\documentclass{amsart}

\usepackage{amsthm,amsfonts,amsmath,amssymb,latexsym,epsfig}

\DeclareMathAlphabet\oldmathcal{OMS}        {cmsy}{b}{n}
\SetMathAlphabet    \oldmathcal{normal}{OMS}{cmsy}{m}{n}
\DeclareMathAlphabet\oldmathbcal{OMS}       {cmsy}{b}{n}
 
\usepackage{eucal}

\usepackage[all]{xy}

\newtheorem{theorem}{Theorem}
\newtheorem{lemma}[theorem]{Lemma}
\newtheorem{proposition}[theorem]{Proposition}

\newtheorem{definition}[theorem]{Definition}
\newtheorem{question}{Question}[section]
\newtheorem{conjecture}[theorem]{Conjecture}
\newenvironment{example}{\medskip \refstepcounter{theorem}
\noindent  {\bf Example \thetheorem}.\rm}{\,}

\renewcommand{\thetheorem}{\thesection.\arabic{theorem}}

\def\d{\partial}                                   
\def\db{\overline{\partial}}
\def\ddb{\partial  \overline{\partial}}  
\def\<{\langle}

\def\>{\rangle}                                    
\def\a{\alpha}

\def\hok{\mbox{}\begin{picture}(10,10)\put(1,0){\line(1,0){7}}
  \put(8,0){\line(0,1){7}}\end{picture}\mbox{}}
\def\BOne{{\mathchoice {\rm 1\mskip-4mu l} {\rm 1\mskip-4mu l}
                          {\rm 1\mskip-4.5mu l} {\rm 1\mskip-5mu l}}}

\def\hook{\mathbin{\hbox to 6pt{%
                 \vrule height0.4pt width5pt depth0pt
                 \kern-.4pt
                 \vrule height6pt width0.4pt depth0pt\hss}}}

\begin{document}

\title{The Sasaki Cone and Extremal Sasakian Metrics}

\author[C.P. Boyer]{Charles P. Boyer}
\author[K. Galicki]{Krzysztof Galicki}
\author[S.R. Simanca]{Santiago R. Simanca}

\thanks{During the preparation of this work, the first two authors
were partially supported by NSF grant DMS-0504367.}

\address{
Department of Mathematics and Statistics, 
University of New Mexico, Albuquerque, N.M. 87131}
\email{cboyer@math.unm.edu, santiago@math.unm.edu}

\begin{abstract}
We study the Sasaki cone of a CR structure of
Sasaki type on a given closed manifold. We introduce an energy functional
over the cone, and use its critical points to single out the strongly
extremal Reeb vectors fields. Should one such vector field be a 
member of the extremal set, the scalar curvature of a Sasaki extremal metric 
representing it would have the smallest $L^2$-norm among all Sasakian 
metrics of fixed volume that can represent vector fields in the cone.
We use links of isolated hypersurface singularities to produce 
examples of manifolds of Sasaki type, many of these in dimension five,  
whose Sasaki cone coincides with the extremal set, and examples where
the extremal set is empty. We end up by proving that a conjecture of 
Orlik concerning the torsion of the homology groups of
these links holds in the five dimensional case.
\end{abstract}

\maketitle

\section{Introduction}
The study of a Sasakian structure goes along with the one dimensional 
foliation on the manifold associated to the Reeb vector field. In the
case where the orbits are all closed, 
the manifold has the structure of an orbifold circle bundle over a 
compact K\"ahler orbifold, which must be algebraic, and which has at most 
cyclic quotient singularities. For a long time, techniques from Riemannian
submersions, suitably extended, were used to find interesting canonical
Sasakian metrics in this orbibundle setting. This led to the common 
belief that the only interesting such metrics occurred precisely in this
setting. Compact Sasakian manifolds with non-closed leaves ---irregular
Sasakian structures--- were 
known to exist, but there was no evidence to suspect they could be had
with Einstein metrics as well. It was thought reasonable
that all Sasaki-Einstein metrics could be understood
well by simply studying the existence of K\"ahler-Einstein metrics on compact
cyclic orbifolds.

The existence of Sasaki-Einstein metrics is of great interest in the physics
of the famous CFT/AdS Duality Conjecture, and the discovery of 
irregular Sasaki-Einstein structures \cite{GMSW04a} was a rather remarkable 
feat, in particular if we take into account their rather 
explicit depiction in coordinates. An attempt to study them
using a treatment inspired in the use of classical Lagrangians, replaced in
this context by a Riemannian functional whose critical points define 
the {\it canonical metrics} we seek, was started in \cite{bgs} using the
squared $L^2$-norm of the scalar curvature as the functional in question, 
with its domain restricted to the space of metrics adapted to 
the geometry of the underlying Sasakian structure under consideration.
This choice of Lagrangian is well known in K\"ahler geometry, which 
was and remains a guiding principle to us, and is quite natural 
since Sasakian geometry sits naturally in between two K\"ahler
geometries, that of the transversals to the one dimensional foliations 
associated to the Reeb vector field, and that of the metric cones inside 
which the Sasaki manifolds sits as a base. Not surprisingly, the 
Sasakian metrics that are critical points of the Lagrangian we use are those
that are transversally extremal in the sense of Calabi, and several important
results known in K\"ahler geometry have now a counterpart in the Sasakian
context \cite{bgs,fow} also.
  
In the theory we developed \cite{bgs}, a given CR structure $({\mathcal D},J)$
of Sasaki type is the analogue in the Sasaki context of a complex structure
on a manifold of K\"ahler type, and the Reeb vector field $\xi$ in a 
Sasakian structure with $({\mathcal D},J)$ as underlying CR structure is the
analogue of a choice of a class in the K\"ahler cone. The 
Sasaki cone $\kappa ({\mathcal D},J)$ of the CR structure arises naturally,
and $\xi\in \kappa ({\mathcal D},J)$ is said to be canonical or extremal if
there exists a Sasakian structure $(\xi,\eta,\Phi, g)$, with underlying 
CR structure $({\mathcal D},J)$, for which the metric $g$ in the structure
is Sasaki extremal. The set of vector fields in the cone 
for which this happens is denoted by ${\mathfrak e}({\mathcal D},J)$.
In this article we present a number of examples where $\kappa ({\mathcal D},J)
= {\mathfrak e}({\mathcal D},J)$, as well as a number of examples where
${\mathfrak e}({\mathcal D},J)$ is the empty set.

In the cases where $\kappa ({\mathcal D},J)={\mathfrak e}({\mathcal D},J)$,
it is natural to ask which is the optimal choice we could make of a Reeb
vector field, be as it may that each one admits a Sasaki extremal 
representative. In effect, this question turns out to be natural 
in general, even in the cases where ${\mathfrak e}({\mathcal D},J)$ is the 
empty set. We briefly discuss this here also, introducing the notion of
strongly extremal Reeb vector field. A detailed elaboration of this
idea will appear elsewhere. 

The specific examples we give are links of isolated hypersurface singularities.
For these manifolds, there is a long standing conjecture of Orlik whose 
validity would describe the torsion of their homology groups. We present 
a proof of this conjecture in dimension five, as obtained recently 
by the second author.

This article developed from partial results obtained by the three
authors while working on an ongoing project. It is presented here
by the first and last author in memory of their friend and collaborator
KRIS GALICKI.

\section{Preliminaries on Sasakian geometry}
\setcounter{theorem}{0}
A contact metric structure $(\xi, \eta, \Phi, g)$ on a manifold $M$ 
is said to be a {\it Sasakian structure} if 
$(\xi, \eta, \Phi)$ is normal. A smooth manifold provided with one such 
structure is said to be a {\it Sasakian manifold}, or a {\it manifold of 
Sasaki type}. We briefly spell out the meaning of this definition in order
to recall the relevant geometric set-up it brings about. 
A detailed discussion of Sasakian and contact geometry can be found 
in \cite{BG05}.

The triple $(\xi, \eta, \Phi)$ defines an almost contact structure
on $M$. That is to say, $\xi$ is a nowhere zero vector field,
$\eta$ is a one form, and $\Phi$ is a tensor of type $(1,1)$, such that
 $$\eta(\xi) = 1 \, , \quad \Phi^2 = -\BOne + \xi \otimes \eta \, .$$
The vector field $\xi$ defines the {\it characteristic foliation} 
${\mathcal F}_{\xi}$ with one-dimensional leaves, and the kernel of $\eta$
defines the codimension one sub-bundle ${\mathcal D}$. We have the
canonical splitting 
 \begin{equation}
 TM = {\mathcal D} \oplus L_{\xi}\, , \label{cs}
 \end{equation}
where $L_\xi$ is the trivial line bundle generated by $\xi.$

The sub-bundle ${\mathcal D}$ inherits an almost complex structure 
$J$ by restriction of $\Phi$, and the dimension of $M$ must
be an odd integer, which we set be $2n+1$ here. 
When forgetting the tensor $\Phi$ and
characteristic foliation, the sub-bundle ${\mathcal D}$ by itself defines
what is called a {\it contact structure} on $M$.

The Riemannian metric $g$ is compatible with the almost contact structure 
$(\xi, \eta, \Phi)$, and so we have have that
$$g(\Phi(X), \Phi(Y))=g(X,Y)-\eta (X)\eta(Y) \, .$$
Thus, $g$ induces an almost Hermitian metric on ${\mathcal D}$. This fact 
makes of $(\xi, \eta, \Phi,g)$ an almost contact metric structure, and
in that case, the canonical decomposition (\ref{cs}) is orthogonal.
The orbits of the field $\xi$ are geodesics, a condition that 
can be re-expressed as $\xi \hok d\eta =0$, and so $d\eta$ is a basic $2$-form.
We have the relation 
 \begin{equation}
 g(\Phi X, Y) = d\eta (X,Y)\, , \label{tkf}
 \end{equation}
which says that $(\xi, \eta, \Phi,g)$ is a contact metric structure.
Notice that the volume element defined by $g$ is given by
 \begin{equation}
 d\mu_g = \frac{1}{n!}\eta \wedge (d\eta)^n  \, . \label{ve}
 \end{equation}

On the cone $C(M)=M\times {\mathbb R}^{+}$ we have the metric 
 $$g_C=dr^2 + r^2 g \, ,$$
and the radial vector field $r\partial_r $ satisfies the relation
 $$\pounds_{r\partial_r} g_C=2g_C \, .$$
We have also an almost complex structure $I$ on $C(M)$ given by
$$I(Y)=\Phi (Y) - \eta(Y) r\partial_r \, , \quad I(r\partial_r ) = \xi \, .$$
The contact metric structure $(\xi, \eta, \Phi)$ is said to be {\it normal} 
if the pair $(C(M),I)$ is a complex manifold. In that case, the induced almost
complex structure $J$ on ${\mathcal D}$ is integrable. Further,
$g_C$ is a K\"ahler metric on $(C(M),I)$.

\subsection{Transverse K\"ahler structure}
For a Sasakian structure $(\xi, \eta, \Phi, g)$, 
the integrability of the almost complex structure $I$ on
the cone $C(M)$ implies that the Reeb vector field $\xi$ leaves 
both, $\eta$ and $\Phi$, invariant. We obtain a codimension one 
integrable strictly pseudo-convex CR structure $({\mathcal D},J)$, where 
${\mathcal D}=\ker \eta$ is the contact bundle and $J=\Phi|_{\mathcal D}$, and 
the restriction of $g$ to ${\mathcal D}$ defines a symmetric form on 
$({\mathcal D},J)$ that we shall refer to as the transverse K\"ahler metric 
$g^T$. By (\ref{tkf}), the K\"ahler form of the transverse K\"ahler metric 
is given by the form $d\eta$. Therefore, the Sasakian metric $g$ is 
determined fully in terms of $(\xi,\eta, \Phi)$ by the expression
 \begin{equation}
 g=d\eta \circ (\BOne \otimes \Phi)+ \eta \otimes \eta \, .\label{me1}
 \end{equation}
The Killing field $\xi$ leaves invariant $\eta$ and $\Phi$, and the 
decomposition (\ref{cs}) is orthogonal. Despite its dependence on the 
other elements of the structure, we insist on explicitly referring to $g$ as 
part of the Sasakian structure $(\xi, \eta, \Phi, g)$. 

We consider the set
\begin{equation}
\label{Sasakianspace}
{\mathcal S}(\xi)=\{\text{Sasakian structure $(\tilde{\xi}, \tilde{\eta},
\tilde{\Phi}, \tilde{g})$} \, | \; \tilde{\xi}=\xi\}\, ,
\end{equation}
and provide it with the $C^\infty$ compact-open topology as sections of
vector bundles.
For any element $(\tilde{\xi}, \tilde{\eta},
\tilde{\Phi}, \tilde{g})$ in this set, the 1-form $\zeta=\tilde{\eta} - \eta$
is basic, and so $[d\tilde{\eta}]_B=[d\eta]_B$. Here, $[\, \cdot \, ]_B$ stands
for a cohomology class in the basic cohomology ring, a ring that is defined by
the restriction $d_B$ of the exterior derivative $d$ to the subcomplex of
basic forms in the de Rham complex of $M$.
Thus, all of the Sasakian structures in ${\mathcal S}(\xi)$ correspond to the
same basic cohomology class. We call
${\mathcal S}(\xi)$ the {\it space of Sasakian structures compatible with
$\xi$}, and say that the Reeb vector field $\xi$ {\it polarizes} the
Sasakian manifold $M$ \cite{bgs}.

Given a Reeb vector field $\xi$, we have its characteristic foliation
${\mathcal F}_\xi$, so we let $\nu({\mathcal F}_\xi)$ be the vector
bundle whose fiber at a point $p\in M$ is the quotient space $T_pM/L_{\xi}$,
and let $\pi_\nu: TM \rightarrow \nu({\mathcal F}_\xi)$ be
the natural projection. The background structure ${\oldmathcal S}=
(\xi,\eta,\Phi,g)$ induces a complex structure $\bar{J}$ on $\nu({\mathcal
F}_\xi)$. This is
defined by $\bar{J}\bar{X}:=\overline{\Phi(X)}$, where $X$ is any
vector field
in $M$ such that $\pi (X)=\bar{X}$. Furthermore, the underlying CR
structure $({\mathcal D},J)$ of ${\oldmathcal S}$ is isomorphic to
$(\nu({\mathcal F}_\xi),\bar{J})$ as a complex vector bundle. For this reason,
we refer to $(\nu({\mathcal F}_\xi),\bar{J})$ as the complex normal bundle
of the Reeb vector field $\xi$, although its identification with
$({\mathcal D},J)$ is not canonical. In this sense, $M$ is polarized by
$(\xi, \bar{J})$.

We define ${\mathcal S}(\xi,\bar{J})$ to be the subset of all
structures $(\tilde{\xi}, \tilde{\eta},
\tilde{\Phi}, \tilde{g})$ in ${\mathcal S}(\xi)$
such that the diagram
\begin{equation}\label{Sasakianspace2eqn}
\begin{array}{ccccc}
TM & \stackrel{\tilde{\Phi}}{\rightarrow} &  TM & \\
  \downarrow \! \! \mbox{{\small $\pi_\nu$}}  & & \downarrow \! \!
\mbox{{\small $\pi_\nu$}} & \\
   \nu({\mathcal F}_\xi) &\stackrel{\bar{J}}{\rightarrow} &\nu(
{\mathcal F}_\xi),&
\end{array}
\end{equation}
commutes. This set consists of elements of ${\mathcal S}(\xi)$
with the same transverse holomorphic structure $\bar{J}$, or with more
precision, the same complex normal bundle $(\nu({\mathcal F}_\xi),\bar{J})$.

\subsection{The Sasaki cone}
If we look at the Sasakian structure $(\xi, \eta, \Phi, g)$ from the 
point of view of CR geometry, its underlying strictly pseudo-convex CR 
structure $({\mathcal D},J)$, with associated contact bundle ${\mathcal D}$, 
has Levi form $d\eta$. 

If $({\mathcal D},J)$ is a strictly pseudo-convex CR structure on $M$ of
codimension one, we say that $({\mathcal D},J)$ is of {\it Sasaki type} if 
there exists a Sasakian structure ${\oldmathcal S}=(\xi,\eta,\Phi,g)$ such 
that 
${\mathcal D}={\rm ker}\, \eta$ and $\Phi|_{\mathcal D}=J$.
We consider the set
\begin{equation}
{\mathcal S}({\mathcal D},J)=\left\{
\begin{array}{c}
 {\oldmathcal S}=(\xi,\eta,\Phi,
g):\; {\oldmathcal S} \; {\rm a\; Sasakian\; structure} \\
({\rm ker}\, \eta, \Phi \mid_{{\rm ker}\, \eta})=
({\mathcal D},J)\end{array} \right\}
\end{equation}
of Sasakian structures {\it with underlying} CR {\it structure} 
$({\mathcal D},J)$.

We denote by 
${\mathfrak c}{\mathfrak o}{\mathfrak n}({\mathcal D})$ the Lie algebra of
infinitesimal contact transformations, and by ${\mathfrak c}{\mathfrak r}
({\mathcal D},J)$ the Lie algebra of the group ${\mathfrak C}{\mathfrak R}
({\mathcal D},J)$ of CR automorphism of
$({\mathcal D},J)$. If 
${\oldmathcal S}=(\xi,\eta,\Phi,g)$ is a contact metric structure whose 
underlying CR structure is $({\mathcal D},J)$, then ${\oldmathcal S}$ is a 
Sasakian structure if, and only if, $\xi \in {\mathfrak c}{\mathfrak r}(
{\mathcal D},J)$.

The vector field $X$ is said to be positive for $({\mathcal D}, J)$ if
$\eta(X)>0$ for any $(\xi,\eta,\Phi,g) \in {\mathcal S}({\mathcal D},J)$. 
We denote by ${\mathfrak c}{\mathfrak r}^{+}({\mathcal D},J)$ the subset of 
all of these positive fields, and consider the projection
\begin{equation}
\begin{array}{ccc}
{\mathcal S}({\mathcal D},J) & \stackrel{\iota}{\rightarrow} &
 {\mathfrak c}{\mathfrak r}^{+}({\mathcal D},J) \\ (\xi,\eta,\Phi,g) & 
\mapsto & \xi
\end{array} \, .\label{ide}
\end{equation}
This mapping identifies naturally
${\mathfrak c}{\mathfrak r}^{+}({\mathcal D},J)$ 
with ${\mathcal S}({\mathcal D},J)$. Furthermore, 
${\mathfrak c}{\mathfrak r}^{+}({\mathcal D},J)$ is an open convex cone 
in ${\mathfrak c}{\mathfrak r}({\mathcal D},J)$ that is invariant under
the adjoint action of the Lie group ${\mathfrak C}{\mathfrak R}({\mathcal D}
,J)$. The identification (\ref{ide}) gives the moduli space 
${\mathcal S}({\mathcal D},J)/{\mathfrak C}{\mathfrak R}({\mathcal D},J)$ of 
all Sasakian structures whose underlying CR structure is $({\mathcal D},J)$.
The conical structure of ${\mathfrak c}{\mathfrak r}^{+}({\mathcal D},J)$ 
justifies the following definition: 

\begin{definition}
Let $({\mathcal D},J)$ be a {\rm CR} structure of Sasaki type on $M$.
The Sasaki cone $\kappa({\mathcal D},J)$ is the moduli space of Sasakian 
structures compatible with $({\mathcal D},J)$,
$$\kappa({\mathcal D},J)={\mathcal S}({\mathcal D},
J)/{\mathfrak C}{\mathfrak R}({\mathcal D},J)\, .$$
\end{definition}

The isotropy subgroup of an element ${\oldmathcal S}\in 
{\mathcal S}({\mathcal D},
J)$ is, by definition, ${\mathfrak A}{\mathfrak u}{\mathfrak t}
({\oldmathcal S})$. It contains a maximal torus $T_k$. In fact, we 
have the following result \cite{bgs}:

\begin{theorem}\label{Sascone}
Let $M$ be a closed manifold of dimension $2n+1$, and let $({\mathcal D},J)$ 
be a {\rm CR} structure of Sasaki type on it. Then the Lie 
algebra ${\mathfrak c}{\mathfrak r}({\mathcal D},J)$ decomposes as 
${\mathfrak c}{\mathfrak r}({\mathcal D},J)= {\mathfrak t}_k+ {\mathfrak p}$,
where ${\mathfrak t}_k$ is the Lie algebra of a maximal torus $T_k$ of 
dimension $k$, $1\leq k\leq n+1$, and ${\mathfrak p}$ is a completely 
reducible $T_k$-module. Furthermore, every $X\in {\mathfrak c}{\mathfrak r}^+
({\mathcal D},J)$ is conjugate to a positive element in the Lie algebra 
${\mathfrak t}_k$. 
\end{theorem}

Theorem \ref{Sascone} and {\rm (\ref{ide})} imply that each orbit can be 
represented by choosing a positive element in the Lie algebra ${\mathfrak t}_k$
of a maximal torus $T_k$. 
So let us fix a maximal torus $T_k$ of a maximal compact subgroup $G$ of 
${\mathfrak C}{\mathfrak R}({\mathcal D},J)$, and let ${\mathcal W}$ denote 
the Weyl group of $G$. If ${\mathfrak t}^{+}_k={\mathfrak t}_k\cap 
{\mathfrak c}{\mathfrak r}^+({\mathcal D},J)$ denotes the subset of positive
 elements in ${\mathfrak t}_k$, we have the identification 
$\kappa({\mathcal D},J)= {\mathfrak t}_k^+/{\mathcal W}$. Each Reeb vector 
field in ${\mathfrak t}_k^+/{\mathcal W}$ corresponds to a unique Sasakian 
structure, so we can view ${\mathfrak t}_k^+/{\mathcal W}$ as a subset of 
${\mathcal S}$, and we have

\begin{lemma}\label{autS}
Let $({\mathcal D},J)$ be a {\rm CR} structure of Sasaki type on $M$, and fix 
a maximal torus $T_k\in {\mathfrak C}{\mathfrak R}({\mathcal D},J)$. Then 
we have that
$$\bigcap_{{\oldmathcal S} \in {\mathfrak t}_k^+/{\mathcal W}} 
{\mathfrak A}{\mathfrak u}{\mathfrak t}_0({\oldmathcal S})= T_k\, ,$$
where ${\mathfrak A}{\mathfrak u}{\mathfrak t}_0({\oldmathcal S})$ 
denotes the identity component of the isotropy group 
${\mathfrak A}{\mathfrak u}{\mathfrak t}({\oldmathcal S})$.
In particular, $T_k$ is contained in the isotropy subgroup of every 
${\oldmathcal S} \in {\mathfrak t}_k^+/{\mathcal W}$.
\end{lemma}

Now the basic Chern class of a Sasakian structure ${\oldmathcal S}
=(\xi,\eta,\Phi, g)$ 
is represented by the Ricci form $\rho^T/2\pi$ of the transverse metric 
$g_T$. Although the notion of basic changes with 
the Reeb vector field, the complex vector bundle ${\mathcal D}$ remains 
fixed. Hence, for any Sasakian structure ${\oldmathcal S}\in 
{\mathcal S}({\mathcal D},J)$, the transverse 
$2$-form $\rho^T/2\pi$ associated to ${\oldmathcal S}$ represents the first 
Chern class  $c_1({\mathcal D})$ of the complex vector bundle 
${\mathcal D}$.

When $k=1$, the maximal torus of ${\mathfrak C}{\mathfrak R}({\mathcal D},J)$ 
is one dimensional. Since the Reeb vector field is central, this implies 
that $\dim {\mathfrak a}{\mathfrak u}{\mathfrak t}({\oldmathcal S}
)= \dim {\mathfrak c}{\mathfrak r}({\mathcal D},J) =1$. Hence, we have that 
${\mathcal S}({\mathcal D},J)={\mathfrak c}{\mathfrak r}^+({\mathcal D},J)=
{\mathfrak t}_1^+={\mathbb R}^+$, and ${\mathcal S}({\mathcal D},J)$ consists 
of the $1$-parameter family of Sasakian structures given by ${\oldmathcal S}_a=
(\xi_a,\eta_a,\Phi_a,g_a)$, where 
\begin{equation}
\xi_a=a^{-1}\xi\, , \quad \eta_a=a\eta\, , \quad \Phi_a=\Phi\, , \quad 
g_a=ag+(a^2-a)\eta\otimes \eta\, , \label{homo}
\end{equation}
and $a\in {\mathbb R}^+$, the $1$-parameter family of transverse homotheties.

In effect, the {\it transverse homotheties} (\ref{homo}) are the only 
deformations 
$(\xi_t,\eta_t,\Phi_t,g_t)$ of a given structure ${\oldmathcal S}=
(\xi,\eta,\Phi,g)$ in 
the Sasakian cone $\kappa({\mathcal D},J)$ where the Reeb vector field varies 
in the 
form $\xi_t=f_t \xi$, $f_t$ a scalar function. For we then have that the 
family of tensors $\Phi_t$ is constant, and since $\pounds_{\xi_t} \Phi_t=0$, 
we see that 
$f_t$ must be annihilated by any section of the sub-bundle ${\mathcal D}$. 
But then (\ref{tkf}) implies that $df_t=(\xi f_t)\eta$, and we conclude that 
the function $f_t$ is constant. Thus, in describing fully the tangent 
space of ${\mathcal S}({\mathcal D},J)$ at ${\oldmathcal S}$, it suffices to
describe only those deformations $(\xi_t,\eta_t,\Phi_t, g_t)$ where
$\dot{\xi}= \partial_t \xi_t \mid_{t=0}$ is $g$-orthogonal to $\xi$. 
These correspond to deformations where the volume of $M$ in the 
metric $g_t$ remains constant in $t$, and are parametrized by 
elements of ${\kappa}({\mathcal D},J)$ that are $g$-orthogonal to
$\xi$.
 
The terminology we use here is chosen to emphasize the analogy that the Sasaki 
cone is to a CR structure of Sasaki type what the K\"ahler cone is to a 
complex manifold of K\"ahler type. Indeed, for any 
point ${\oldmathcal S}=(\xi,\eta,\Phi,g)$ in $\kappa({\mathcal D},J)$, the 
complex normal bundle $(\nu({\mathfrak F}_{\xi}),\bar{J})$ is isomorphic to 
$({\mathcal D},J)$.
In this sense, the complex structure $\bar{J}$ is fixed with the fixing of  
$({\mathcal D},J)$, the Reeb vector field $\xi$ polarizes the manifold, and 
the Sasaki cone $\kappa({\mathcal D},J)$ represents the set of all possible 
polarizations. We observe though that when we fix 
$\xi \in \kappa({\mathcal D},J)$, the 
underlying CR 
structure of elements in ${\mathcal S}(\xi,\bar{J})$ might change, even though 
their normal bundles are all isomorphic to $({\mathcal D},J)$.

Let $\xi \in \kappa({\mathcal D},J)$. We may ask if there are 
canonical representatives of ${\mathcal S}(\xi,\bar{J})$. 
We can answer this question using a variational principle, as 
done in \cite{bgs}. 
For let ${\mathfrak M}(\xi,\bar{J})$ be the set of all 
compatible Sasakian metrics arising from structures in 
${\mathcal S}(\xi, \bar{J})$, and consider the
functional
\begin{equation}
\begin{array}{ccl}
{\mathfrak M}(\xi,\bar{J})  & \rightarrow & {\mathbb R}\, , \\
g & \mapsto & {\displaystyle \int _M s_g ^2 d{\mu}_g }\, .
\end{array}
\label{var}
\end{equation}
Its critical points define the canonical representatives that we seek.
Furthermore, by the identification of 
$\kappa({\mathcal D},J)$ with ${\mathfrak t}_k^+$, we may then
single out the set of Reeb vector fields $\xi$ in ${\mathfrak t}_k^+$ 
for which the functional (\ref{var}) admits critical points at all.

\begin{definition}
We say that $(\xi,\eta,\Phi,g)\in {\mathcal S}(\xi,\bar{J})$
is an extremal Sasakian structure if 
$g$ is a critical point of {\rm (\ref{var})}. A vector field
$\xi \in \kappa({\mathcal D},J)$ is said to be extremal if
there exists an extremal Sasakian structure in ${\mathcal S}(\xi,\bar{J})$.
We denote by ${\mathfrak e}({\mathcal D},J)$ the set of all 
extremal elements of the Sasaki cone, and refer to it as the 
extremal Sasaki set of the {\rm CR} structure $({\mathcal D},J)$.
\end{definition}

Notice that $(\xi,\eta,\Phi,g)\in 
{\mathcal S}(\xi,\bar{J})$ is extremal if, and only if, the 
transversal metric $g^T$ is extremal in the K\"ahler sense \cite{bgs}.
We also have that ${\mathfrak e}({\mathcal D},J)$ is an open subset of
$\kappa({\mathcal D},J)$ \cite{bgs}, a Sasakian version of the
the openness theorem in K\"ahler geometry \cite{ls,si}.

\section{The energy functional in the Sasaki cone}
\setcounter{theorem}{0}
By Lemma \ref{autS}, given any Sasakian structure $(\xi ,\eta,\Phi,g)$ in
${\mathcal S}({\mathcal D},J)$, the space of $g^T$-Killing potentials
 associated to a maximal torus
has dimension $k-1$. We denote by ${\mathcal H}_{\xi }=
{\mathcal H}_{\xi,g}({\mathcal D},J)$ the vector subspace of $C^{\infty}(M)$ 
that 
they and the constant functions span. This is the space of functions whose
transverse gradient is holomorphic. 

Given a Sasakian metric $g$ such that 
$(\xi ,\eta,\Phi,g) \in {\mathcal S}({\mathcal D},J)$, we define the map
\begin{equation}
\begin{array}{ccl}
C^{\infty}(M) & \rightarrow & {\mathcal H}_{\xi } \\
f & \mapsto & \pi_{g} f 
\end{array} 
\end{equation}
to be the $L^2$-orthogonal projection with ${\mathcal H}_{\xi }$ as its
range. Then the metric $g$ in ${\mathfrak M}(\xi,\bar{J})$ is extremal if,
and only if, the scalar curvature $s_g=s_g^T-2n$ is equal to its 
projection onto ${\mathcal H}_{\xi }$ \cite{si2}. Or said differently, if, 
and only if, $s_g$ is an affine function of the space of Killing potentials.
The functional (\ref{var}) then admits the lower bound
$${\mathfrak M}(\xi,\bar{J}) \ni g \mapsto \int s_g^2 d\mu_g \geq
\int (\pi_{g}s_g)^2 d\mu_g \, ,$$
where the right side only depends upon the polarization $(\xi,{\mathcal D},J)$,
and it is only reached if the metric $g$ is extremal. This 
reproduces in this context a situation analyzed in the K\"ahler case already
\cite{si3,si4}.

We define the energy of a Reeb vector field in the Sasakian cone to be the
functional given by this optimal lower bound:
\begin{equation}\label{en}
\begin{array}{ccl}
\kappa({\mathcal D},J) & \stackrel{E}{\rightarrow} & {\mathbb R} \\
\xi & \mapsto & {\displaystyle \int (\pi_{g}s_g)^2 d\mu_g }\, ,$$ 
\end{array} 
\end{equation}
 
\begin{definition}
A Reeb vector field $\xi \in \kappa ({\mathcal D},J)$ is said to be
strongly extremal if it is a critical point of the functional {\rm (\ref{en})}
over the space of Sasakian structures in $\kappa ({\mathcal D},J)$ that fix 
the volume of $M$. The Sasakian structure
$(\xi ,\eta,\Phi,g)$ is said to be a strongly extremal representative
of $\kappa({\mathcal D},J)$ if $\xi $ is strongly extremal, and if $g
\in {\mathcal S}(\xi, \bar{J})$ is an extremal representative of the
polarized manifold $(M,\xi, J)$. 
\end{definition}

\section{The Euler-Lagrange equation for strongly extremal 
Reeb vector fields}
\setcounter{theorem}{0}
We now discuss some of the variational results that play a r\^ole in the
derivation of the Euler-Lagrange equation of a strongly extremal Reeb vector
field. Complete details of this derivation will appear elsewhere.

Let us consider a path $\xi_t \in \kappa({\mathcal D},J)$, and a
corresponding path of Sasakian structures $(\xi_t ,\eta_t,\Phi_t,g_t)$.
The path of contact forms must be of the type 
\begin{equation}
\eta_t = \frac{1}{\eta(\xi_t)}\eta + \frac{1}{2}d^c \varphi_t \, , \label{cf}
\end{equation}
where $\varphi_t$ is a $\xi_t$-basic function. For convenience, we drop the
subscript when we consider the value of quantities at $t=0$. 
Then, we have that $\eta(\xi)=1$, $\varphi$ is constant, and 
$\dot{\varphi}$ a $\xi$-basic function.  

\begin{proposition}\label{pr1}
Let $(\xi_t ,\eta_t,\Phi_t,g_t)$ be a path of Sasakian structures that 
starts at $(\xi ,\eta,\Phi,g)$, such that 
$\xi_t \in \kappa({\mathcal D},J)$ and $\eta_t$ is of the form
{\rm (\ref{cf})}. Then we have that
$$d\mu_t = 1-t\left((n+1)\eta(\dot{\xi})+\frac{1}{2}\Delta_B \dot{\varphi} 
\right) d\mu +O(t^2) \, ,$$
and
$$\frac{d}{dt}\mu_{g_t}(M)=-(n+1) \int (\eta(\dot{\xi}))\, d\mu \, .$$
Volume preserving Sasakian deformations are given by variations $\xi_t$ of the
Reeb vector field whose infinitesimal deformation is globally orthogonal to
it, that is to say, such that $\eta(\dot{\xi})$ is orthogonal to the constants.
Finally, 
$$D^{2}\mu (\dot{\xi},\dot{\beta})=(n+1)(n+2)\int (\eta(\dot{\xi}))\,
(\eta(\dot{\beta}))\, d\mu \, .$$
\end{proposition}

{\it Proof}. By differentiation of (\ref{ve}), we obtain the variational 
expression for the volume form. The expression for the infinitesimal
variation of the volume follows by Stokes' theorem, and the one for the 
Hessian by an iteration of this argument. \qed

\begin{proposition} \label{pr2}
Let $(\xi_t ,\eta_t,\Phi_t,g_t)$ be a path of Sasakian structures that 
starts at $(\xi ,\eta,\Phi,g)$, such that 
$\xi_t \in \kappa({\mathcal D},J)$ and $\eta_t$ is of the form
{\rm (\ref{cf})}. 
Then we have the expansion
$$
s_t = s^T-2n -t\left(n\Delta_B (\eta(\dot{\xi}))-s^T \eta(\dot{\xi})+
\frac{1}{2}\Delta_B^2 \dot{\varphi}+2(\rho^T,
i\ddb \dot{\varphi})\right) +O(t^2) \, ,
$$
for the scalar curvature of $g_t$. 
\end{proposition}

{\it Proof}. We have the relation $s_t=s^T_t -2n$, where $s^T_t$ is the
scalar curvature of the transverse K\"ahler metric $d\eta_t$. Since
the latter is 
$$s^T_t=-2(d\eta_t, i\ddb \log{\det{(d\eta_t)}})_{d\eta_t}\, ,$$
taking the $t$-derivative, and evaluating at $t=0$, we obtain the desired
result. Just notice that the $d\eta$-trace of $d\dot{\eta}$ is equal
to the $d\eta$-trace of $-\eta(\dot{\xi})d\eta+i\ddb \dot{\varphi}$,
as the trace of $d(\eta(\dot{\xi} )) \wedge \eta$ is zero. \qed 

We may now sketch the computation of the variation of the
functional (\ref{en}) along the curve $\xi_t$. If the subscript $t$ denotes
quantities associated to the Sasakian structure $(\xi_t ,\eta_t,\Phi_t,g_t)$,
we have that
$$\frac{d}{dt}E(\xi_t)  =  2\int (\pi_t s_t )
(\frac{d}{dt}(\pi_t s_t ))d\mu_t + \int (\pi_t s_t)^2  \frac{d}{dt} d\mu_t
\, ,$$
and by the results of the previous propositions, the value at $t=0$ is
given by
\begin{equation}
\frac{d}{dt}E(\xi_t)\mid_{t=0}  =  2\int (\pi s )\left[
(-n\Delta_B +s^T)\eta(\dot{\xi}) +\dot{\pi} s \right] d\mu
 -(n+1)\int (\pi s)^2 \eta(\dot{\xi})d\mu\, . \label{t1}
\end{equation}

We may compute the variation of $\pi$, but need less than that for our
purposes. Indeed, if $G_g$ is the Green's operator of the metric $g$, we have
that
$$
\int (\pi s) \dot{\pi}s\, d\mu=\int (s^T-\pi s^T)2i G_g
\overline{\d}_g^{*}(\d^{\#}\pi s\hok d\dot{\eta})\, ,
$$
where given a basic function $f$, $\partial^{\#}f$ is the operator obtained 
by raising the indices of the $(0,1)$ form $\db f$ using the transversal 
K\"ahler metric $d\eta$, and since
$$i\d^{\#}_g\pi s \hok d\dot{\eta}=-\db (\pi s \, \eta(\dot{\xi}) ) \, ,$$ 
we obtain that 
\begin{equation}
\int (\pi s) \dot{\pi}s\, d\mu=\int (s^T-\pi s^T)2i G_g
\overline{\d}_g^{*}(\d^{\#}\pi s\hok d\dot{\eta})=-
\int (s^T-\pi s^T) \pi s \, \eta(\dot{\xi}) d\mu \, , \label{t2}
\end{equation}
where the last equality follows because $s^T-\pi s^T$ is orthogonal to
the constants. 

Combining these results, we obtain the following:

\begin{theorem} \label{count}
Let $(\xi_t ,\eta_t,\Phi_t,g_t)$ be a path of volume preserving 
Sasakian structures that starts at $(\xi ,\eta,\Phi,g)$, such that 
$\xi_t \in \kappa({\mathcal D},J)$, and $\eta_t$ is of the form
{\rm (\ref{cf})} with $\varphi_t$ a constant. Then we have that
\begin{equation}
\frac{d}{dt}E(\xi_t)\mid_{t=0}  =  2\int (\pi s )\left[
(-n\Delta_B +\pi s^T)\eta(\dot{\xi}) \right] d\mu
 -(n+1)\int (\pi s)^2 \eta(\dot{\xi})d\mu\, . \label{t3}
\end{equation}
The Reeb vector field $\xi \in \kappa({\mathcal D},J)$ is a critical
point of {\rm (\ref{en})} over the space of deformations with $g_t$ of
fixed volume if, and only if, 
\begin{equation}
\pi_g \left[ -2n \Delta_B \pi_g s_g -(n-1)(\pi_g s_g)^2\right] +
4n \pi_g s_g =\lambda \, ,
\label{el}
\end{equation}
for $\lambda$ a constant. \qed
\end{theorem}

\subsection{Some examples of critical Reeb vector fields}
We may define a character that obstructs the existence of extremal Sasakian 
structures $(\xi, \eta, \Phi, g) \in {\mathcal S}(\xi, \bar{J})$ where the
scalar curvature $s_g$ of $g$ is constant \cite{bgs}. This is analogous 
to the Futaki invariant \cite{fu,ca2} in the K\"ahler category.
 
Indeed, for any transversally holomorphic vector
field $\d^{\#}_g f$ with Killing potential $f$, the Sasaki-Futaki 
invariant is given by
\begin{equation}
{\mathfrak F}_{\xi}(\d_g^{\#}f)=- \int f(s_g-s_{0})d\mu_g \, , \label{fi}
\end{equation}
where $g$ is any Sasakian metric in ${\mathfrak M}(\xi,\bar{J})$, $s_g$ is
its scalar curvature, $s_0$ is the projection of $s_g$ onto the constants.
This expression defines uniquely a character on the entire algebra of 
transversally holomorphic vector fields, and if $(\xi, \eta, \Phi, g)$
is an extremal Sasakian structure in ${\mathcal S}(\xi, \bar{J})$, 
the scalar curvature $s_g$ is constant if, and only
if, ${\mathfrak F}_{\xi} \equiv 0$.

\begin{proposition}
Suppose that $(\xi, \eta, \Phi, g) \in {\mathcal S}(\xi, \bar{J})$ is such 
that $\pi_g s_g$ is a constant. Then, for any other Sasakian structure
$(\xi,\tilde{\eta},\tilde{\Phi},\tilde{g})$ in ${\mathcal S}(\xi, \bar{J})$, 
we have that $\pi_{\tilde{g}} s_{\tilde{g}}$ is constant also, and equal to 
$\pi_g s_g$.
\end{proposition}

{\it Proof}. The proof follows by the same reasoning used to prove
an analogous statement in the K\"ahler context (cf. \cite{lusa1}, 
Proposition 5). \qed

\begin{theorem}
Let $\xi \in \kappa( {\mathcal D},J)$ be such that the Sasaki-Futaki
invariant is identically zero, ${\mathfrak F}_{\xi}( \, \cdot \, )=0$. 
Then $\xi$ is a critical point of  {\rm (\ref{en})}. 
\end{theorem}

{\it Proof}. If ${\mathfrak F}_{\xi}( \, \cdot \, )=0$, then $\pi_g s_g$ is 
constant, and $\xi$ satisfies (\ref{el}). \qed

\begin{example}
Under the assumption that the cone $C(M)=M\times {\mathbb R}^{+}$ is a 
Calabi-Yau manifold, the Hilbert action over Sasakian metrics has been  
extensively analyzed \cite{MaSpYau06}. When considered as a function of 
Reeb vector fields of charge $n+1$, this functional is a multiple of
the volume\footnotemark, and its variation over such a domain is a constant
multiple 
of the Sasaki-Futaki invariant. Thus, a critical Reeb vector field $\xi$ has
vanishing Sasaki-Futaki invariant, and if $(\xi,\eta, \Phi, g)$ is a Sasakian
structure corresponding to it, the function $\pi_g s_g$ is 
necessarily constant. 
By the convexity of the volume functional (see Proposition \ref{pr1}), 
the critical Reeb vector field of this functional is unique.
We thus see that the critical Reeb field singled out by the Hilbert action,
on this type of special manifolds,  must be a critical point of 
our functional (\ref{en}) also, as was to be expected. \qed
\end{example}

\footnotetext{In earnest, this is done in \cite{MaSpYau06} assuming that 
$M$ is quasi-regular, but the result holds in further generality. Notice also
that in the statement for the {\it charge} of the Reeb vector field, we have 
made adjustments between the convention we use for the dimension of $M$,
and the one in the said reference.} 

\section{Sasakian geometry of links of isolated hypersurface 
singularities}
\setcounter{theorem}{0}
Many of the results we discuss from here on involve the natural Sasakian 
geometry occurring on links of isolated hypersurface singularities. In the
spirit of self-containment, we provide a brief review of the geometry of 
links in arbitrary dimension, which we shall use mostly in dimension $5$
later on.

Let us recall (cf. Chapters 4 and 9 of \cite{BG05}) 
that a polynomial $f\in {\mathbb C}[z_0,\ldots,z_n]$ is said to be a 
{\it weighted homogeneous polynomial} of {\it degree} $d$ and 
{\it weight} ${\bf w}=(w_0,\ldots,w_n)$ if for any $\lambda \in 
{\mathbb C}^{\times}$
$$
f(\lambda^{w_0} z_0,\ldots,\lambda^{w_n}
z_n)=\lambda^df(z_0,\ldots,z_n)\, .
$$
We are interested in those weighted homogeneous polynomials $f$ whose 
zero locus in ${\mathbb C}^{n+1}$ has only an isolated singularity at the 
origin. We define the {\it link} $L_f({\bf w},d)$ as 
$f^{-1}(0)\cap {\mathbb S}^{2n+1}$, where ${\mathbb S}^{2n+1}$ is the 
$(2n+1)$-sphere in Euclidean space. By the Milnor fibration theorem 
\cite{Mil68}, $L_f({\bf w},d)$ is a closed $(n-2)$-connected manifold that 
bounds a parallelizable manifold with the homotopy type of a bouquet of 
$n$-spheres. Furthermore, $L_f({\bf w},d)$ admits a 
Sasaki-Seifert structure 
${\oldmathcal S}=(\xi_{\bf w},\eta_{\bf w},\Phi_{\bf w},g_{\bf w})$ in a 
natural way 
\cite{Tak,BG01b}. This structure is quasi-regular, and its 
bundle $({\mathcal D},J)$ has $c_1({\mathcal D})=0$. This latter
property implies that 
\begin{equation}
c_1({\mathcal F}_{\xi_{\bf w}})=a[d\eta_{\bf w}]_B \label{sig}
\end{equation}
for some constant $a$, where ${\mathcal F}_{\xi_{\bf w}}$ is the 
characteristic foliation. The sign of $a$ determines the negative, null, 
and positive cases that we shall refer to below. 

The Sasaki-Futaki invariant (\ref{fi}) is the obstruction for the scalar
curvature of the metric $g$ in an extremal Sasakian structure 
${\oldmathcal S}=(\xi, \eta, \Phi, g)$ to be constant \cite{bgs}. 
One particular case where this obstruction vanishes is when there are no 
transversally holomorphic vector fields other than those generated by 
$\xi$ itself. In that case, the Sasaki cone of the corresponding CR structure 
turns out to be of dimension one. If the Sasakian structure is quasi-regular,
the transversally holomorphic vector fields are the 
lifts of holomorphic vector fields on the projective algebraic variety 
${\mathcal Z}=M/{\mathbb S}^1_{\xi}$, where ${\mathbb S}^1_{\xi}$ is the 
circle generated by the Reeb vector field.

Let us recall that a Sasakian metric $g$ is said to be $\eta$-Einstein if
$${\rm Ric}_g=\lambda g+ \nu \eta\otimes \eta $$ for some 
constants $\lambda$ and $\nu$ that satisfy the relation
$\lambda+\nu =2n$. They yield particular examples of Sasakian metrics of
constant scalar curvature, a function which in these cases is given in terms 
of $\lambda$ by the expression $s_g=2n(1+\lambda)$.
 
The only transversally holomorphic vector fields of a negative or a
null link of an isolated hypersurface singularity are those
generated by the Reeb vector field. These links have a one-dimensional
Sasaki cone, and by the transverse Aubin-Yau theorem 
\cite{BGM06}, any point in this cone can be represented by a Sasakian
structure whose metric is $\eta$-Einstein (see examples of such links 
in \cite{BoGa05a,BG05}). We thus obtain the following 

\begin{theorem}
Let $L_f({\bf w},d)$ be either a negative or null link of an isolated 
hypersurface singularity with underlying {\rm CR}
 structure $({\mathcal D},J)$.
Then its Sasaki cone $\kappa({\mathcal D},J)$ is one-dimensional,
and it coincides with the extremal Sasaki set ${\mathfrak e}({\mathcal D},J)$. 
\end{theorem}

The case of positive links is more complicated. This can be deduced already 
from the fact that if we were to have a Sasakian structure 
${\oldmathcal S}=(\xi, \eta, \Phi, g)$ such that 
$c_1({\mathcal F}_{\xi})=a[d\eta]_B$ with $a>0$ and where $g$ is 
$\eta$-Einstein, then there would exist a transverse homothety (\ref{homo}) 
that would yield a Sasaki-Einstein metric of positive scalar curvature, and 
generally speaking, these types of metrics are difficult to be had.

In order to understand this case better, we begin by stating the following
result that applies to links in general, and which follows by Corollary 5.3 
of \cite{bgs}:

\begin{proposition}\label{linkexthyp}
Let $L_f({\bf w},d)$ be the link of an isolated hypersurface singularity 
with its natural Sasakian structure ${\oldmathcal S}= 
(\xi_{\bf w},\eta_{\bf w},\Phi_{\bf w},g_{\bf w})$. If 
$g\in {\mathcal S}(\xi_{\bf w},\bar{J})$ is an extremal Sasakian metric with 
constant scalar curvature, then $g$ is $\eta$-Einstein.
\end{proposition}

By Lemma 5.5.3 of  \cite{BG05}, if the 
weight vector ${\bf w}=(w_0,\ldots,w_n)$ of a link $L_f({\bf w},d)$ is such 
that $2w_i< d$ for all but at most one of the indices, then the 
Lie algebra ${\mathfrak a}{\mathfrak u}{\mathfrak t}({\oldmathcal S})$
is generated by the Reeb vector field $\xi_{\bf w}$.
We have the following

\begin{theorem}\label{nosym}
Let $L_f({\bf w},d)$ be the link of an isolated hypersurface singularity 
with $2w_i< d$ for all but at most one of the indices, 
and let $({\mathcal D},J)$ by its underlying {\rm CR} structure.
Then the dimension of the Sasaki cone $\kappa({\mathcal D},J)$ is one.
\end{theorem}

Even if we assume that for a given positive link $L_f({\bf w},d)$ the 
hypothesis of Theorem \ref{nosym} apply, and so its Sasaki-Futaki invariant 
vanishes identically, there could be other 
obstructions preventing the existence of extremal Sasakian metrics 
of constant scalar curvature on it. Indeed, as was observed in 
\cite{GMSY06}, classical estimates of Bishop and Lichnerowicz may obstruct 
the existence of Sasaki-Einstein metrics, an observation that has
produced a very effective tool to rule out the existence of these metrics
in various
cases. Employing the Lichnerowicz estimates \cite{GMSY06} in combination with
the discussion above and some results in \cite{bgs}, we obtain the following:

\begin{proposition}\label{Lich.Ob.Link.Ineq} Let $L_f({\bf w},d)$ be a 
link of an isolated hypersurface singularity with its natural Sasakian 
structure, and let $I:=|{\bf w}|-d=(\sum_j w_j) -d$ be its index. If 
$$I> n \, \min_{i}{\{w_i\}}\, ,$$
then $L_f({\bf w};d)$ cannot admit any Sasaki-Einstein metric. In that case,
${\mathcal S}(\xi_{\bf w},\bar{J})$ does not admit extremal 
representatives.
\end{proposition}

Combining Theorem \ref{nosym} and Proposition \ref{Lich.Ob.Link.Ineq}, we 
may obtain examples of positive links whose Sasaki cones are one dimensional,
and whose extremal sets are empty.

\begin{theorem}
Let $L_f({\bf w},d)$ be a positive link of an isolated hypersurface singularity
whose index $I$ satisfies the relation 
$I> n \, \min_{i}{\{w_i\}}$, and such that 
$2w_i< d$ for all but at most one of the indices. 
If $({\mathcal D},J)$ is its underlying {\rm CR} structure, we have
that the dimension of the Sasaki cone $\kappa({\mathcal D},J)$ is one,
and that the extremal set ${\mathfrak e}({\mathcal D},J)$ is empty.
\end{theorem}

Besides the conditions given above on the weight vector of a link,
there are other interesting algebraic conditions that suffice to prove
the existence of extremal Sasakian metrics on it. These arise when we view
the link as a Seifert fibration over a base $X$ with orbifold 
singularities. If the K\"ahler class of the projective 
algebraic variety $X$ is a primitive integral class in the second 
orbifold cohomology group, and the constant $a$ in (\ref{sig}) is positive, 
they serve to measure the singularity of the pair 
$(X,K^{-1}+\Delta)$, where $K^{-1}$ is an anti-canonical divisor, and 
$\Delta$ is
a branch divisor, and are known as 
{\it Kawamata log terminal} (Klt) conditions (cf. Chapter 5 of \cite{BG05}). 

A weighted homogeneous polynomial of the form 
\begin{equation}
f=z_0^{a_0}+\cdots +z_n^{a_n}\, , \quad a_i \geq 2\, ,
\label{bp}
\end{equation}
is called a {\it Brieskorn-Pham} (BP) polynomial. In this case the 
exponents $a_i$, the weights $w_i$, and degree are related by $d=a_iw_i$ 
for each $i=0,\ldots,n.$ We change slightly the 
notation in this case, and denote by $L_f({\bf a})$, 
${\bf a}=(a_0, \ldots , a_n)$, the link that a BP polynomial 
$f$ defines. These are special but quite important examples of links. Their
Klt conditions were described in \cite{BGK05}. The 
base of a BP link 
$L_f({\bf a})$ admits a positive K\"ahler-Einstein orbifold metric if
\begin{equation}\label{bgkest}
 1<\sum_{i=0}^n\frac{1}{a_i}<
1+\frac{n}{n-1}\min_{i,j}\Bigl\{\frac{1}{a_i},
\frac{1}{b_ib_j}\Bigr\}\,.
\end{equation}
where $b_j=\gcd(a_j,C^j)$ and $C^j={\rm lcm}\{ a_i : \, i\neq j\}$.
This condition leads to the finding of a a rather large 
number of examples of Sasaki-Einstein metrics on homotopy spheres
\cite{BGK05,BGKT05} and rational homology spheres \cite{BG06b,Kol05b}.

In the special case when the integers $(a_0,\ldots,a_n)$ are pairwise 
relatively prime, Ghigi and Koll\'ar \cite{GhKo05} obtained a sharp 
estimate. In this case, the BP link is always a homotopy sphere, and if
we combine the now sharp Klt estimate with Proposition \ref{Lich.Ob.Link.Ineq}
above, we see that, when the $a_i$s are pairwise relatively prime, a BP
link $L_f({\bf a})$ admits an extremal 
Sasakian metric if, and only if,  
$$\sum_{i=0}^n\frac{1}{a_i}<
1+n~\min_{i}\Bigl\{\frac{1}{a_i}\Bigr\}\, .$$

Other applications of the Klt estimate (\ref{bgkest}) exist 
\cite{BoGa05a,BG05}. If $f$ is neither a BP polynomial nor a perturbation 
thereof, alternative estimates must be developed. Further details can be
found in \cite{BG05}.

Notice that in all the link examples we have discussed above, 
information about the topology of the underlying manifold is absent.
We end this section addressing that issue in part.

Topological information about links was first given by Milnor and Orlik in 
\cite{MiOr70} through the study of the Alexander polynomial. They computed
 the Betti numbers of the manifold underlying the link of any weighted 
homogeneous polynomial of an isolated hypersurface singularity, and in the 
case of a rational homology sphere, 
the order of the relevant homology groups. 
A bit later, Orlik \cite{Or72} postulated a combinatorial 
conjecture for computing the torsion, an algorithm that we describe next
(see Chapter 9 of \cite{BG05} for more detail). 

Given a link $L_f({\bf w}, d)$, we define its {\it fractional weights}
to be 
\begin{equation}\label{frct.1}
\Bigl(\frac{d}{w_0},\cdots, \frac{d}{w_n}\Bigr)\equiv
\Bigl(\frac{u_0}{v_0},\cdots, \frac{u_n}{v_n}\Bigr),
\end{equation}
where
\begin{equation}\label{uv.wd.conv}
u_i=\frac{d}{{\rm gcd}(d,w_i)},\qquad v_i=\frac{w_i}{{\rm
gcd}(d,w_i)}.
\end{equation}
We denote by $({\bf u}, {\bf v})$ the tuple $(u_0,\ldots,u_n,
v_0,\ldots,v_n)$. By \eqref{uv.wd.conv}, we may 
go between $({\bf w},d)$ and $({\bf u}, {\bf v})$. 
We will sometimes write $L_f({\bf u}, {\bf v})$ for 
$L_f({\bf w}, d)$.
   
\begin{definition}\label{orl.alg}{\bf (Orlik's algorithm)}
Let $\{i_1,\ldots,i_s\} \subset\{0,1,\ldots,n\}$ be an ordered set of 
$s$ indices, that is to say, $i_1<i_2<\cdots<i_s$. Let us denote by $I$ its 
power set {\rm (}consisting of all of the $2^s$ subsets of the set{\rm )}, 
and by $J$ the set of all proper subsets. 
Given a $(2n+2)$-tuple $({\bf u}, {\bf v})=(u_0,\ldots,u_n,
v_0,\ldots,v_n)$ of integers, we define inductively a set of $2^s$ positive 
integers, one for each ordered element of $I$, as follows: 
$$c_{\emptyset}=\gcd{( u_0,\ldots, u_n )}\, ,$$
and if $\{i_1,\ldots,i_s\}\in I$ is ordered, 
then
\begin{equation}\label{Orlikc's}
c_{i_1,\ldots,i_s}=\frac{\gcd{(u_0,\ldots,\hat{u}_{i_1},\ldots,\hat{u}_{i_s},
\ldots,u_n)}}{\prod_J c_{j_1,\ldots j_t}}\, .
\end{equation}
Similarly, we also define a set of $2^s$ real numbers by
$$k_\emptyset=\epsilon_{n+1}\, ,$$
and
\begin{equation}
k_{i_1,\ldots,i_s}= \epsilon_{n-s+1}\sum_{J}(-1)^{s-t}
\frac{u_{j_1}\cdots u_{j_t}}{v_{j_1}\cdots v_{j_t}{\rm
lcm}( u_{j_1},\ldots,u_{j_t})}\, ,
\end{equation}
where
$$\epsilon_{n-s+1}=\left\{ \begin{array}{cl}
0 & \text{if $n-s+1$ is even,} \\
1 & \text{if $n-s+1$ is odd,}
\end{array}
\right.
$$
respectively. 
Finally, for any $j$ such that $1\leq j\leq r=[ \max{\{k_{i_1,\ldots,i_s}\}}]$,
where $[x]$ is the greatest integer less or equal than $x$, we set
\begin{equation}\label{torsionorders}
d_j=\prod_{ k_{i_1,\ldots,i_s}\geq j}c_{i_1,\ldots,i_s}\, .
\end{equation}\qed
\end{definition}

{\it Orlik's torsion conjecture} \cite{Or72} is stated in terms of the
integers computed by this algorithm:

\begin{conjecture}\label{orlikconjtor} For a link $L_f({\bf u},{\bf v})$ with
fractional weights $({\bf u},{\bf v})$, we have  
$$H_{n-1}(L_f({\bf u},{\bf v}),{\mathbb Z})_{\rm tor}=
{\mathbb Z}/d_1 \oplus \cdots \oplus {\mathbb Z}/d_r \, .$$
\end{conjecture}

\noindent This conjecture is known to hold in certain special cases 
\cite{Ran75,OrRa77a}.
\begin{proposition}\label{Orlikconjknown}
{\rm Conjecture \ref{orlikconjtor}} holds in the following  cases:
\begin{enumerate}
\item In dimension $3$, that is to say, when $n=2.$ 
\item For Brieskorn-Pham polynomials {\rm (\ref{bp})}.
\item For $f({\bf z})=z_0^{a_0}+z_0z_1^{a_1}+z_1z_2^{a_2}+\cdots
+z_{n-1}z_n^{a_n}$. 
\end{enumerate}
\end{proposition}

\noindent It is also known to hold for certain complete intersections given by 
generalized Brieskorn polynomials \cite{Ran75}. 

Before his tragic accident, by using Koll\'ar's Theorem 
\ref{simpconn.over.ratsurf.thm} below, Kris proved that Conjecture 
\ref{orlikconjtor} holds in dimension 5 also.
We present his argument in the appendix.

\section{Extremal Sasakian metrics in dimension five}
\setcounter{theorem}{0}
We start by reviewing briefly some general facts about Sasakian geometry 
on simply connected 5-manifolds (cf. Chapter 10 of \cite{BG05}).

\subsection{Sasakian geometry in dimension five}
Sasakian geometry in dimension five is large enough to be interesting
while remaining manageable, at least in the case of closed simply connected 
manifolds where we can use the Smale-Barden classification. 

Up to diffeomorphism, Smale \cite{Sm62} classified all closed simply 
connected $5$-manifolds that admit a spin structure, and showed they  
must be of the form
\begin{equation}\label{smaleman}
M=kM_\infty\# M_{m_1}\#\cdots \#M_{m_n}
\end{equation} 
where $M_{\infty}= {\mathbb S}^2\times {\mathbb S}^3$, 
$kM_{\infty}$ is the $k$-fold connected sum of $M_{\infty}$,   
$m_i$ is a positive integer with $m_i$ dividing $m_{i+1}$ and $m_1\geq 1$, 
and where $M_{m}$ is ${\mathbb S}^5$ if $m=1$, or a 
$5$-manifold such that $H_2(M_m,{\mathbb Z})={\mathbb Z}/m \oplus 
{\mathbb Z}/m$, otherwise. 
The integer $k$ in this expression can take on the values 
$0,1,\ldots $, with $k=0$ corresponding to the case where there is
no $M_{\infty}$ factor at all. It will be convenient to
set the convention $0M_{\infty}={\mathbb S}^5$ below, which is consistent
with the fact that the sphere is the neutral element for the connected
sum operation. The $m_i$s can range in $1,2,\ldots $.

Barden \cite{Bar65} extended Smale's  classification to include the
non-spin case also, where we must add the non-trivial 
${\mathbb S}^3$-bundle over ${\mathbb S}^2$, denoted by $X_\infty$, 
and certain rational 
homology spheres $X_j$ parametrized by $j=-1,1,2,\ldots$ that  
do not admit contact structures, and can be safely ignored from consideration
when studying the Sasakian case. Notice that 
$X_\infty\# X_\infty=2M_\infty$, so in considering simply connected Sasakian
$5$-manifolds that are non-spin, it suffices to take a connected sum 
of the manifold (\ref{smaleman}) with one copy of $X_\infty$, or 
equivalently, replace at most one copy of $M_\infty$ with $X_\infty$.

What Smale-Barden manifolds can admit Sasakian structures? A partial answer 
to this question has been given recently by Koll\'ar \cite{Kol06a}. In order
to present this, it is convenient to describe the torsion subgroup of 
$H_2(M,{\mathbb Z})$ for the manifold $M$ in terms of elementary divisors 
rather than the invariant factors in (\ref{smaleman}). 
The group $H_2(M,{\mathbb Z})$ can be written as a direct sum of cyclic 
groups of prime power order
\begin{equation}\label{BaSm-dec}
H_2(M,{\mathbb Z})={\mathbb Z}^k\oplus \bigoplus_{p,i}
\bigl({\mathbb Z}/{p^i}\bigr)^{c(p^i)}\, ,
\end{equation}
where $k=b_2(M)$ and $c(p^i)=c(p^i,M)$. These non-negative integers
are determined by $H_2(M,{\mathbb Z})$, but the subgroups
${\mathbb Z}/{p^i}\subset H_2(M,{\mathbb Z})$ are not unique. 
We can choose the decomposition (\ref{BaSm-dec}) such that the second
Stiefel-Whitney class map
$$
w_2:H_2(M,{\mathbb Z})\to {\mathbb Z}/2
$$ 
is zero on all but one of the summands  ${\mathbb Z}/{2^j}$. 
If we now assume that $M$ carries a Sasakian structure, 
by Rukimbira's approximation theorem \cite{Ruk95a}, $M$ admits a 
Seifert fibered structure (in general, non-unique) with an associated 
Sasakian structure, or a {\it Sasaki-Seifert structure} as referred to
in \cite{BG05}. Koll\'ar \cite{Kol06a} proved that the existence of a 
Seifert fibered structure on $M$ imposes constraints on the invariants 
$k$ and $c(p^i)$ of (\ref{BaSm-dec}) above. Namely, for the primes 
that appear in the expression, the cardinality of the set 
$\{i :\; c(p^i)>0\}$ must be less or equal than $k+1$.
For example, both $M_2$ and $M_4$ admit Sasakian 
structures but $M_2\#M_4$ does not. Other obstructions involving the torsion 
subgroup of  $H_2(M,{\mathbb Z})$ apply as well \cite{Kol06a, BG05}.

The case of positive Sasakian structures is particularly attractive 
since they provide examples of Riemannian metrics with positive Ricci 
curvature. Koll\'ar \cite{Kol05b} (see also \cite{Kol06} in these
Proceedings) has shown that positivity greatly restricts 
the allowable torsion groups. He proved that if a closed simply 
connected 5-manifold $M$ admits a positive Sasakian structure, then the 
torsion subgroup of $H_2(M,{\mathbb Z})$ must be one of the following:
\begin{enumerate}
\item $({\mathbb Z}/m)^2$, $m\in {\mathbb Z}^+$, 
\item $({\mathbb Z}/5)^4$, 
\item $({\mathbb Z}/4)^4$, 
\item $({\mathbb Z}/3)^4$, $({\mathbb Z}/3)^6$, or $({\mathbb Z}/3)^8$,
\item $({\mathbb Z}/2)^{2n}$, $n\in {\mathbb Z}^+$.
\end{enumerate} 
Conversely, for each finite group $G$ in the list above, there is
a closed simply connected 5-manifold $M$ with $H_2(M,{\mathbb Z})_{\rm tor}=G$,
and remarkably enough, all of these manifolds can be realized as the links 
of isolated hypersurface singularities.

Let us recall that a quasi-regular Sasakian structure 
can be viewed as a Seifert fibered structure. In dimension five,
Koll\'ar \cite{Kol05b} has shown how the branch divisors $D_i$ of the 
orbifold base determine the torsion in $H_2(M,{\mathbb Z})$. 
For simply connected 5-manifolds, this result is as follows:

\begin{theorem} \label{simpconn.over.ratsurf.thm}
Let $M^5$ be a compact simply connected 5-manifold with a quasi-regular 
Sasakian structure ${\oldmathcal S}$, let $({\mathcal Z},
\sum_i (1-\tfrac1{m_i})D_i)$ denote 
the corresponding projective algebraic orbifold base, with branch divisors 
$D_i$ and ramification index $m_i$, and let $k$ be the second Betti number 
of $M^5$. Then the integral cohomology groups of $M^5$ are as follows: 
\medskip
\begin{center}
\begin{tabular}{|c|c|c|c|c|c|c|}\hline \hline
i & {\rm 0} & {\rm 1} & {\rm 2} & {\rm 3} & {\rm 4} & {\rm 5} \\
\hline \hline
$H^i(M^5,{\mathbb Z})$& ${\mathbb Z}$ & $0$& ${\mathbb Z}^{k}$ &
${\mathbb Z}^{k}\oplus\sum_i
({\mathbb Z}/{m_i})^{2g(D_i)}$ & $0$ &${\mathbb Z}$ \\
\hline
\end{tabular}
\end{center}
\medskip

\noindent where $g(D)$ is the genus of the Riemann surface $D$
{\rm (}notice that ${\mathbb Z}/1$ is the trivial group{\rm )}.
\end{theorem}

\subsection{Existence of Extremal Metrics}
Here we reproduce a table \cite{BG05} that lists all simply connected 
spin 5-manifolds that can admit a Sasaki-Einstein metric (and which are, 
therefore, Sasaki extremal). In the first column we list the 
type of manifold in terms of Smale's description (\ref{smaleman}), 
while in the second column we indicate restrictions on these under which a 
Sasaki-Einstein structure is known to exist. Any Smale manifold that is not 
listed 
here cannot admit a Sasaki-Einstein metric, but could, in principle, admit an 
extremal metric.
We list the manifolds but not the number of deformation classes of positive
 Sasakian 
structures that may occur, a number that varies depending upon the manifold 
in question. For example, there are infinitely many such deformation classes 
on the 5-sphere ${\mathbb S}^5$ and on $k({\mathbb S}^2\times {\mathbb S}^3)$,
as well as on some rational homology spheres. On the other hand, there is a 
unique deformation class of positive Sasakian structures on the rational 
homology spheres $2M_5$ and $4M_3$ \cite{Kol05b,Kol06}, and they both 
admit extremal metrics.

\medskip
\begin{center}
\begin{tabular}{|c|c|}\hline
\hline 
Manifold $M$ & S-E\\
\hline\hline $kM_\infty$, $k\geq0$ & any $k$ \\
\hline
$8M_\infty\#M_m$, $m>2 $ & $m>4$ \\
\hline
$7M_\infty\# M_m$, $m>2$ & $m>2$\\
\hline
$6M_\infty\# M_m$, $m>2$ & $m>2$\\
\hline
$5M_\infty\# M_m$, $m>2$ & $m>11$\\
\hline
$4M_\infty\# M_m$, $m>2$& $m>4$\\
\hline
$3M_\infty\# M_m$, $m>2$ & $m=7,9$  or $m>10$ \\
\hline
$2M_\infty\# M_m$, $m>2$ & $m>11$ \\
\hline
$M_\infty\# M_m$, $m>2$ & $m>11$\\
\hline
$M_m$, $m>2$ & $m>2$ \\
\hline $2M_5$, $2M_4$, $4M_3$, $M_\infty\#2M_4$ & yes \\
\hline
$kM_\infty\# 2M_3 $ & $k=0$\\
\hline
$kM_\infty\# 3M_3$, $k\geq0$ & $k=0$\\
 \hline
$kM_\infty\# nM_2$, $k\geq 0$, $n>0$ & $(k,n)=(0,1)$ or $(1,n)$, $n>0$\\
\hline $kM_\infty\# M_m$, $k>8$, $2<m<12$ & \\
 \hline
\end{tabular}
\vspace{3mm}\\

\parbox{4.00in}{\small Table 1. Simply connected spin 5-manifolds admitting
 Sasaki-Einstein metrics. The right column indicates the restriction that 
ensures that the manifold
on the left carries a Sasaki-Einstein metric.}
 
\end{center}
\medskip

\subsection{Brieskorn-Pham links with no extremal Sasakian metrics}
We now provide a table with examples of Brieskorn-Pham links 
$L_f({\bf a})$ whose extremal sets are all empty. These are obtained by 
using Proposition \ref{Lich.Ob.Link.Ineq}, the required 
estimate on the index being easily checked.
We list the link together with the range for its parameters, 
the underlying manifold, and the dimension of the Sasaki cone 
${\kappa}$ of the corresponding CR structure. 
There are infinitely many such links (or Sasaki-Seifert structures) 
on each of the manifolds listed. 
\medskip
 
\begin{center}
\begin{tabular}{|c|c|c|}\hline
\hline
Manifold $M$ & Link $L_f({\bf a})$ & ${\rm Dim}\, {\kappa}$ \\
\hline\hline
$M_1$ & $L(2,3,3,6l+1)$, $l\geq 2$ & $1$\\
\hline
$M_1$ & $L(2,3,3,6l+5)$, $l\geq 2$ & $1$\\
\hline
$M_1$ & $L(2,2,2,2l+1)$, $l\geq 2$ & $2$\\
\hline
$M_\infty$ & $L(2,2,2,2l)$, $l\geq 3$ & $2$\\
\hline
$2M_\infty$ & $L(2,3,3,2(3l+1))$, $l\geq 2$ & $1$\\
\hline
$2M_\infty$ & $L(2,3,3,2(3l+2))$, $l\geq 2$ & $1$\\
\hline
$4M_\infty$ & $L(2,3,3,6l)$, $l\geq 3$ & $1$\\
\hline
$6M_\infty$ & $L(2,3,4,12l)$, $l\geq 3$ & $1$\\
\hline
$8M_\infty$ & $L(2,3,5,30l)$, $l\geq 3$ & $1$\\
\hline
$(k-1)M_\infty$ & $L(2,2,k,2lk)$, $l\geq 2$, $k\geq 3$ & $2$\\
\hline
$4M_2$ & $L(2,3,5,15(2l+1))$, $l\geq 2$ & $1$\\
\hline
$2M_\infty\# M_2$ & $L(2,3,4,6(2l+1))$, $l\geq 2$ & $1$\\
\hline
\end{tabular}
\vspace{3mm}\\

\parbox{4.00in}{\small Table 2. Some examples of 5-dimensional Brieskorn-Pham 
links whose associated space of Sasakian structures, ${\mathcal S}(\xi)$,
does not have extremal Sasakian metrics. 
The last column
lists the dimension of their Sasaki cones.}

\end{center}
\medskip

More examples can be had by considering general weighted homogeneous 
polynomials $f$ instead. In those cases, we go 
back to the notation $L_f({\bf w},d)$ for the link, ${\bf w}$ its weight 
vector, and $d$ its degree. Their systematic study 
should be possible by using the classification of normal forms of 
Yau and Yu \cite{YaYu05}. Here we content ourselves by giving a few examples.
In all of these cases, the Sasaki cone is one dimensional.
\medskip

\begin{center}
\begin{tabular}{|c|c|c|}\hline
\hline
Manifold $M$ & weight vector ${\bf w}$ & $d$ \\
\hline\hline
$3M_\infty\# M_2$ & $(2,2(4l+3),6l+5,2(6l+5))$, $l\geq 2$ & $4(6l+5)$ \\
\hline
$7M_\infty$ & $(1,4l+1,3l+1,2(3l+1))$, $l\geq 2$ & $4(3l+1)$ \\
\hline
$M_4$ & $(4,3(2l+1),4(2l+1),4(3l+1))$, $l\geq 5$ & $12(2l+1)$ \\
\hline
\end{tabular}
\vspace{3mm}\\

\parbox{4.00in}{\small Table 3. 5-dimensional links of 
isolated hypersurface singularities that do not carry Sasakian extremal
structure.}
\end{center}
\medskip

In the first row of this table the existence of positive Sasakian structures 
on the Smale-Barden manifold $3M_\infty\# M_2$ is reported here for the 
first time. As indicated in the table, for $l\geq 2$ these links do not 
carry extremal Sasakian metrics. It is unknown if an extremal
Sasakian metric exists when $l=1$, and generally speaking, it is not known 
whether an extremal Sasakian metric exists at all on $3M_\infty\# M_2$. 
A similar situation occurs for the
link $L_f((4,5,12,20), 40)$, which gives a positive Sasakian structure on 
the manifold $3M_\infty\# M_4$, and which is not yet known to carry 
extremal Sasakian metrics. Thus, it is natural to pose the following:

\begin{question}
Are there positive Sasakian manifolds that admit no extremal Sasakian metrics?
\end{question}

\subsection{Toric Sasakian $5$-manifolds}
The study of contact toric $(2n+1)$-manifolds goes back to \cite{BM93}, 
where it was shown that toric contact structures split into two 
types, those where the torus action is free, and those where it is not,
with the free action case easily described, and the
non-free case described via a Delzant type theorem when $n\geq 2$.
A somewhat refined classification was obtained later \cite{Ler02a},
and the non-free case (again, for $n\geq 2$) was characterized in
\cite{BG00b} by the condition that the Reeb vector field belongs 
to the Lie algebra of the $(n+1)$-torus acting on the manifold. 
These toric actions were called actions of {\it Reeb type}, and it was proved 
that every toric contact structure of Reeb type admits a 
compatible Sasakian structure. Here we content ourselves with describing the 
Sasaki cone for a simple but interesting example. For more on toric Sasakian 
geometry, we refer the reader to \cite{BGO06,fow,CFO07}.

\begin{example}
The Wang-Ziller manifold $M^{p_1,p_2}_{k_1,k_2}$ \cite{WaZi90}
is defined to be the total space of the ${\mathbb S}^1$-bundle over 
${\mathbb C}{\mathbb P}^{p_1}\times {\mathbb C}{\mathbb P}^{p_2}$
whose first Chern class is $k_1\a_1+k_2\a_2$, $\a_i$ being
the positive generator of $H^2({\mathbb C}{\mathbb P}^{p_i}, {\mathbb Z})$
with $k_1, k_2 \in {\mathbb Z}^+$. These manifolds all admit homogeneous 
Einstein metrics \cite{WaZi90}, and also admit homogeneous Sasakian structures
\cite{BG00a}. However, they are Sasaki-Einstein only if 
$k_1\a_1+k_2\a_2$ is proportional to the first Chern class of 
${\mathbb C}{\mathbb P}^{p_1}\times {\mathbb C}{\mathbb P}^{p_2}$, that is to
say, only when $k_1\a_1+k_2\a_2$ and $(p_1+1)\a_1+ (p_2+1)\a_2$ are 
proportional.

 We are interested in $M^{p,q}_{k,l}$ as a toric Sasakian manifold, 
and treat the five dimensional case $M^{1,1}_{k_1,k_2}$ only.
We shall also assume that $\gcd(k_1,k_2) =1$, for then 
$M^{1,1}_{k_1,k_2}$ is simply connected, and 
diffeomorphic to ${\mathbb S}^2\times {\mathbb S}^3$. 
If $k_1$ and $k_2$ are not 
relatively prime, the analysis of $M^{1,1}_{k_1,k_2}$ can be obtained easily 
from the simply connected case by taking an appropriate cyclic quotient. We 
shall also assume that if $(k_1,k_2)\neq (1,1)$, then $k_1<k_2$, for 
otherwise we can simply interchange the two ${\mathbb C}{\mathbb P}^{1}$s.

We wish to find the Sasaki cone for 
$M^{1,1}_{k_1,k_2}$. For that, we identify this manifold with the quotient 
${\mathbb S}^3\times {\mathbb S}^3/{\mathbb S}^1(k_1,k_2)$,
where the ${\mathbb S}^1$ action on 
${\mathbb S}^3\times {\mathbb S}^3$ is given by
$$({\bf x},{\bf y})\mapsto (e^{ik_2\theta}{\bf x},e^{-ik_1\theta}{\bf y})\, ,$$
and points in ${\mathbb S}^3$ are viewed as points in ${\mathbb C}^2$.
Alternatively, $M^{1,1}_{k_1,k_2}$ can be obtained
by a Sasakian reduction \cite{GrOr01} of the standard 
Sasakian structure ${\oldmathcal S}_0$ on ${\mathbb S}^7$ by the 
${\mathbb S}^1(k_1,k_2)$ action, the moment map   
$\mu: {\mathbb S}^7 \rightarrow {\mathbb R}$ being given by 
$$\mu =k_2(|z_0|^2+|z_1|^2) -k_1(|z_2|^2+|z_3|^2)\, .$$
The zero level set is the product ${\mathbb S}^3(R_1)\times {\mathbb S}^3(R_2)$
of spheres whose radii are given by
\begin{equation}\label{radii}
R_i^2=\frac{k_i}{k_1+k_2}\, .
\end{equation}
The Sasakian structure on $M^{1,1}_{k_1,k_2}$ induced by this reduction 
is denoted by ${\oldmathcal S}_{k_1,k_2}=(\xi,\eta,\Phi,g)$, where for ease of 
notation, we refrain from writing the tensor fields with the subscripts
also.

Let us consider coordinates $(z_0, z_1, z_2, z_3)$ in 
${\mathbb C}^4$, where $z_j=x_j+iy_j$. Then the vector fields
$$H_i=y_i\frac{\partial}{\partial x_i}-x_i\frac{\partial}{\partial y_i}\, , 
\; 0\leq i \leq 3\, ,$$
restrict to vector fields on ${\mathbb S}^7$ that span the Lie algebra 
${\mathfrak t}_4$ of the maximal torus $T^4\subset {\mathbb S}{\mathbb O}(8)$.
The base space ${\mathbb C}{\mathbb P}^{1}\times 
{\mathbb C}{\mathbb P}^{1}$ of $M_{k_1,k_2}^{1,1}$ is obtained
from the zero level set ${\mathbb S}^3(R_1)\times {\mathbb S}^3(R_2)$ 
as the quotient by the free $T^2$-action given by 
$({\bf x},{\bf y})\mapsto (e^{i\theta_1}{\bf x},e^{i\theta_2}{\bf y})$.
Let ${\tilde X}$ denote the vector field on ${\mathbb C}{\mathbb P}^{1}$ 
generating the ${\mathbb S}^1$-action defined in homogeneous coordinates by 
$[{\zeta}_0,{\zeta}_1]\mapsto [e^{i\phi}{\zeta}_0,
e^{-i\phi}{\zeta}_1]$, and denote by ${\tilde X}_1$ and ${\tilde X}_2$ the 
vector field ${\tilde X}$ on the two copies of ${\mathbb C}{\mathbb P}^{1}$.
We let $X_1$ and $X_2$ denote the lifts of these vector fields
to ${\mathfrak a}{\mathfrak u}{\mathfrak t}({\oldmathcal S}_{k_1,k_2})$.

The Reeb vector field $\xi$ of ${\oldmathcal S}_{k_1,k_2}$ is that induced by 
the restriction of $H_0+H_1+H_2+H_3$ to the zero level set 
${\mathbb S}^3(R_1)\times {\mathbb S}^3(R_2)$ of the moment map $\mu$
(see Theorems 8.5.2 and 8.5.3 in \cite{BG05}). The 
Lie algebra ${\mathfrak t}_3(k_1,k_2)$ of the  3-torus of 
${\oldmathcal S}_{k_1,k_2}$ is spanned by the Reeb vector 
field $\xi$, and the two vector fields $X_1$ and $X_2$ described above.

We have:
\begin{lemma}
The Sasaki cone ${\kappa}({\mathcal D}_{k_1,k_2},J)$ for the Wang-Ziller 
manifold $M^{1,1}_{k_1,k_2}$ is given by the set of all vector fields
$$Z_{l_0,l_1,l_2}=l_0\xi +l_1X_1+l_2X_2 \, , \; (l_0,l_1,l_2)\in 
{\mathbb R}^3\, ,$$
subject to the conditions 
$$l_0>0,\qquad |k_1l_1+k_2l_2|<(k_1+k_2)l_0.$$
\end{lemma}

{\it Proof}. The Sasaki cone $\kappa ({\mathcal D}_{k_1,k_2},J)$ is equal 
to the set
of $\{Z\in {\mathfrak t}_3 \, | \; \eta(Z)>0\}$ of positive vector fields
(see Theorem \ref{Sascone}).
The vector fields $X_1$ and $X_2$ are induced on the quotient 
$M^{1,1}_{k_1,k_2}\approx ({\mathbb S}^3(R_1)\times {\mathbb S}^3(R_2))/
{\mathbb S}^1(k_1,k_2)$ by the restrictions of
$H_0-H_1$ and $H_2-H_3$ to ${\mathbb S}^3(R_1)\times {\mathbb S}^3(R_2)$.
So we have
$$\eta(Z_{l_0,l_1,l_2})=\eta(l_0\xi+l_1X_1+l_2X_2)= l_0+l_1(|z_0|^2-|z_1|^2) 
+l_2(|z_2|^2-|z_3|^2)\, ,$$
with $l_0>0$. The desired result now follows if we use the expression 
(\ref{radii}) for the radii. \qed

For $(l_0,l_1,l_2)=(1,0,0)$, we obtain the original K\"ahler structure on 
the base ${\mathbb C}{\mathbb P}^1\times {\mathbb C}{\mathbb P}^1$,  
with K\"ahler metric $k_1\omega_1+k_2\omega_2$ where $\omega_i$ is the
standard Fubini-Study K\"ahler form on the corresponding 
${\mathbb C}{\mathbb P}^1$ factor. The scalar curvature of this metric 
is constant, and therefore, the Reeb vector field $\xi$ on $M^{1,1}_{k_1,k_2}$ 
is strongly extremal. Moreover, by the openness theorem \cite{bgs}, there
exist an open neighborhood of $\xi$ in 
$\kappa({\mathcal D}_{k_1,k_2},J)$ that is
entirely contained in the extremal set 
${\mathfrak e}({\mathcal D}_{k_1,k_2},J)$.
It would be of interest to determine how big a neighborhood this can be,
and how it compares to the entire Sasaki cone.
Although each $M^{1,1}_{k_1,k_2}$ carries a Sasakian structure 
${\oldmathcal S}_{k_1,k_2}$ with an extremal Sasakian metric of 
constant scalar curvature, 
the only one with a Sasaki-Einstein metric is $M^{1,1}_{1,1}$.
This follows from the fact that $c_1({\mathcal D}_{k_1,k_2})=2(k_2-k_1)x$,
where $x$ is the positive generator of $H^2({\mathbb S}^2\times {\mathbb S}^3,
{\mathbb Z})$. 
\end{example}

\section{Appendix}
We begin describing a graphical procedure that generalizes the now  
well-known Brieskorn graphs (see Theorem 10.3.5 and Remark 10.3.1 of 
\cite{BG05} for a detailed discussion).

To a given link $L_f({\bf w},d)=L_f({\bf u},{\bf v})$ defined by a weighted 
homogeneous polynomial $f$, of weight vector ${\bf w}$ and fractional weight 
vector $({\bf u},{\bf v})$, we associate a graph as follows:

\begin{definition}\label{B.gra.rat.5}
Consider the following Rational Brieskorn Polytope
$$
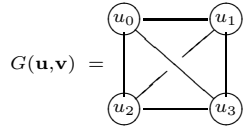
\xymatrix{ *+[o][F-]{\scriptstyle{u_0}}\ar@{-}[d]_{G({\bf u},{\bf v})\; \; =\;
\; }\ar@{-}[dr]\hole\ar@{-}[r]
             & *+[o][F-]{\scriptstyle{u_1}} \ar@{-}[d]\ar@{-}[dl]|\hole \\
       *+[o][F-]{\scriptstyle{u_2}} & *+[o][F-]{\scriptstyle{u_3}}\ar@{-}[l] }
$$
where we think of $G({\bf u},{\bf v})$ as a tetrahedron, and
\begin{enumerate}
\item we label the vertices with pairs $(u_i,\alpha_i)$, $i=0,1,2,3$,
where $\alpha_i=1-1/v_j-1/v_k-1/v_l$ for a set of distinct indices 
$\{i,j,k,l\}$,
\item we label the edges with the 
numbers $\alpha_{ij}=\gcd(u_i,u_j)/v_jv_j$,
\item we label the faces with numbers $\alpha_{ijk}=\gcd(u_i,u_j,u_k)/
v_iv_jv_k$,
\item we label the interior of the tetrahedron with the 
pair 
$$(t,\tau)\equiv \Bigl(
\frac{u_0u_1u_2u_3}{v_0v_1v_2v_3 {\rm
lcm}(u_1,u_1,u_2,u_3)}, -1+\frac{1}{v_0}+ \frac{1}{v_1}+\frac{1}{v_2}+
\frac{1}{v_3}\Bigr)\, .
$$
\noindent Furthermore,
\item we define the reduced indices $m_i$ to be the
factor of $u_i$ that is not a factor of any of the $u_j$s associated to
the remaining $3$ vertices,
\item we define numbers $g_i$ so that $2g_i+\alpha_i$ is computed 
from the face opposite to the vertex $(u_i,\alpha_i)$ as the product of 
the edge numbers divided by the face number minus the sum of the edge
numbers,
\item we define $\kappa$ so that $\kappa+\tau-t$ is the sum of the
six edge numbers minus the sum of four products of
edge numbers divided by the face numbers, one such term for
each face.
\end{enumerate}
\end{definition}

Using this associated graph, we can reformulate 
Orlik's Conjecture \ref{orlikconjtor} in dimension 
five as follows:

\begin{conjecture}\label{orlik5}
The second homology group over ${\mathbb Z}$ of a five dimensional link 
$L_f({\bf u},{\bf v})$ is given by
$$H_2(L({\bf a}),{\mathbb Z})={\mathbb Z}^\kappa
\oplus({\mathbb Z}/m_{0})^{2g_0}\oplus({\mathbb Z}/m_1)^{2g_1}
\oplus({\mathbb Z}/m_2)^{2g_2}\oplus({\mathbb Z}/m_3)^{2g_3}\, ,$$
where $\kappa , m_i, g_i$ are the integers given in {\rm Definition 
\ref{B.gra.rat.5}}.
\end{conjecture}

That $\kappa$ is the second Betti number of the link was known long ago
\cite{MiOr70}. There remains to describe the torsion.

Before going any further in this direction, we need to reformulate 
Koll\'ar's Theorem \ref{simpconn.over.ratsurf.thm} for 5 dimensional 
links $L_f({\bf w},d)$. First, 
let us recall the following two simple conditions of
quasi-smoothness for curves \cite{Fle00} and surfaces
\cite{JoKo01a}:

\begin{proposition}\label{QS.curves} {\rm \cite{Fle00}} 
Let $C_d\subset {\mathbb C}{\mathbb P}^2(w_0,w_1,w_2)$ be a curve of 
degree $d$ defined by a homogeneous polynomial $f$. Assume that $d>a_i$,
$a_i$ integer, $i=0,1,2$. Then $C_d$ 
is quasi-smooth if for all $i\in\{0,1,2\}$ the following
conditions are satisfied:
\begin{enumerate}
\item the polynomial $f$ has a monomial $z_i^{a_i}z_j$ for some $j$ {\rm (}here
we allow $i=j${\rm )} of  degree $d$
\item the polynomial $f$ has a monomial of degree $d$ which does not 
involve $z_i$.
\end{enumerate}
\end{proposition}

\begin{proposition}\label{QS.for.surfaces}{\rm \cite{Fle00,JoKo01a}}
A general hypersurface $X_d \subset {\mathbb C}{\mathbb P}(w_0,w_1,w_2,w_3)$ 
defined by a homogeneous polynomial $f$ is quasi-smooth if,
 and only if, all of the following three conditions hold 
\begin{enumerate}
\item For each $i=0,\ldots,3$ there is an index $j$ such that $f$ has a 
monomial
$z_i^{m_i}z_j\in H^0({\mathbb C}{\mathbb P}({\bf w}),{\mathcal O}(d))$.
Here $j=i$ is possible.
\item If $\gcd(w_i,w_j)>1$ then $f$ has a monomial
$z_i^{b_i}z_j^{b_j}\in H^0({\mathbb C}{\mathbb P}({\bf w}),{\mathcal O}(d))$.
\item For every $i,j$ either $f$ has a monomial
$z_i^{b_i}z_j^{b_j}\in H^0({\mathbb C}{\mathbb P}({\bf w}),{\mathcal O}(d))$,
 or $f$ has monomials $z_i^{c_i}z_j^{c_j}z_k$ and $z_i^{d_i}z_j^{d_j}z_l\in
H^0({\mathbb C}{\mathbb P}({\bf w}),{\mathcal O}(d))$ with 
$\{k,l\}\neq \{i,j\}$.
\end{enumerate}
\end{proposition}

We are then ready for the following proposition.

\begin{proposition}\label{prop.kol.links}
Let $L_f({\bf w},d)$ be a link defined by a quasi-smooth weighted
homogeneous polynomial $f=f(z_0,z_1,z_2,z_3)$, with weights
${\bf w}=(w_0,w_1,w_2,w_3)$, and of degree $d$. Let
$L_f \rightarrow X_f \subset {\mathbb C}{\mathbb P}^3(w_0,w_1,w_2,w_3)$ 
be the associated Seifert fibration. Consider the sets
$$D_i=X_f\cap\{z_i=0\}\subset {\mathbb C}{\mathbb P}^2
(w_0,\ldots,\hat{w_i},\ldots,w_3)\, .$$
Note that $D_i$ need not be an orbifold, so pick $i\in
A\subset\Omega=\{0,1,2,3\}$ such that $D_i$ is a quasi-smooth
{\rm (}orbifold{\rm )} Riemann surface. Then $D_i$ can be singular only
if $\gcd{( w_0,\ldots,\hat{w}_i,\dots,w_3)}=1$, and we have that 
$$H_2(L_f,{\mathbb Z})={\mathbb Z}^{\kappa}\oplus\sum_{i}
({\mathbb Z}/m_i)^{2g(D_i)}\, ,$$ 
where $\kappa=b_2(L_f)=b_2(X_f)-1$, and the
ramification indices are given by 
$m_i=\gcd{\{w_0,\ldots,\hat{w_i},\ldots,w_3\}}$.
In particular, let us set $d_i=d/m_i$, and let us normalize the weight vector
by $\tilde{{\bf w}}^{i}=
\frac{1}{m_i}(w_0,\ldots,\hat{w_i},\ldots,w_3)$. 
Then $D_i \subset {\mathbb C}{\mathbb P}^2(\tilde{\bf w}^i)$ is an
orbifold Riemann surface of degree $d_i$ embedded in the weighted
projective $2$-space with weights $\tilde{\bf w}^{i}$, and its 
genus $g(D_i)$ is given by the formula \cite{OrWa71a,Or72b}
\begin{equation}\label{genusformula}
g(\Sigma_{(\tilde{\bf w}, d)})=\frac{1}{2}\Bigl(\frac{d^2}{
\tilde{w_0}\tilde{w_1}\tilde{w_2}}-d\sum_{i<j}
\frac{\gcd(\tilde{w_i},\tilde{w_j})}{\tilde{w_i}\tilde{w_j}}+
\sum_i\frac{\gcd(d,\tilde{w_i})}{\tilde{w_i}}-1\Bigr)\, .
\end{equation}
\end{proposition}

{\it Proof}. We only need to show that if
$D_i$ is singular, then  we must have $m_i=1$ so that $D_i\subset
X_f^{reg}$ lies inside the smooth part of $X_f$. This will follow from the
quasi-smoothness condition on $f$.

Suppose then that $m_0=\gcd(w_1,w_2,w_3)>1$, and that $D_0$ is
singular. We shall derive a contradiction. 

Since $m_0>1$, we have $\gcd(w_i,w_j)>1$ for all 
$1\leq i<j\leq3$. Since $X_f \subset {\mathbb C}{\mathbb P}^3
(w_0,w_1,w_2,w_3)$ is quasi-smooth, $f({\bf z})$ must contain the 
following monomial terms (other than the terms involving $z_0$):
$$z_1^{a_1}z_j, \quad z_2^{a_2}z_k, \quad z_3^{a_3}z_l,\qquad
j,k,l\in\{1,2,3\}\, ,$$ and
$$z_1^{\alpha}z_2^{\beta},\quad z_2^\gamma z_3^\delta, \quad z_3^\tau
z_1^\rho\, .$$ 
Note that, for example, $z_1^{a_1}z_0$ is not possible
as this would mean that 
$$d=w_0+w_1a_1\, ,$$
and $m_0$ divides $d$ and $w_1$ but not $w_0$. Hence, the
polynomial $g(z_1,z_2,z_3)=f(0,z_1,z_2,z_3)$ will contain these
six monomial terms. Note that the degree of $g$ is 
$\tilde{d}=d/m_0$, and it has weights $\tilde{w}_i=w_i/m_0$.
By Proposition \ref{QS.curves}, it is quasi-smooth. This contradicts
the fact that $D_0$ is assumed to be singular. \qed

\subsection{Proof of Orlik's Conjecture in Dimension $5$}
In order to prove Conjecture \ref{orlik5}, we only need to show that 
the $(\kappa, m_i,g_i)$s of Definition \ref{B.gra.rat.5} are the same as 
the ones in Proposition \ref{prop.kol.links}. We need to pass between 
the $({\bf u},{\bf v})$-data and the $({\bf w},d)$-data. 

It suffices to prove the needed statement for $i=0$. We first observe 
that $(m_0,k_0)=(c_{123},k_{123}).$ In
terms of the $({\bf u},{\bf v})$-data, we have
$$m_0=\frac{u_0{\rm gcd}(u_0,u_1,u_2){\rm gcd}(u_0,u_1,u_2){\rm
gcd}(u_0,u_2,u_3){\rm gcd}(u_0,u_1,u_3)}{{\rm
gcd}(u_0,u_1,u_2,u_3){\rm gcd}(u_0,u_1){\rm gcd}(u_0,u_2) {\rm
gcd}(u_0,u_3)},$$
and
$$\begin{array}{rcl}
2g_0 & = & {\displaystyle 
-1+\Bigl(\frac{1}{v_1}+ \frac{1}{v_2}+\frac{1}{v_3}\Bigr)-
\Bigl(\frac{u_1u_2}{v_1v_2{\rm lcm}(u_1,u_2)}+
\frac{u_1u_3}{v_1v_3{\rm lcm}(u_1,u_3)} +\frac{u_2u_3}{v_2v_3{\rm
lcm}(u_2,u_3)}\Bigr) } \vspace{1mm} \\
& & +{\displaystyle  \frac{u_1u_2u_3}{v_1v_2v_3{\rm
lcm}(u_1,u_2,u_3)}} \, .
\end{array}
$$
We recall now that 
$$u_i=\frac{d}{\gcd({d,w_i)}}\, , \qquad 
v_i=\frac{w_i}{{\gcd({d,w_i) }}}\, ,
\qquad \frac{u_i}{v_i}=\frac{d}{w_i}\, .$$

Let us begin with $2g_0$.  The first term is no problem as
$$\frac{1}{v_i}=\frac{\gcd(d,w_i)}{w_i}\, ,$$
as needed. 

Since $f$ is quasi-smooth, $\gcd(w_i,w_j)$ must divide the degree
for all $i<j$. By this, we have the following
$$\begin{array}{rcl}
{\displaystyle \frac{u_iu_j}{v_iv_j{\rm lcm}({u_i,u_j)} }} & = &
{\displaystyle  \frac{\gcd({u_i,u_j)}}{v_iv_j}=\frac{\gcd({u_i,u_j)}
\gcd({d,w_j)} \gcd({d,w_i)}}{w_iw_j} } \vspace{1mm} \\
 & = & {\displaystyle  \frac{\gcd\bigl({ u_i (\gcd({ d,w_j)} 
\gcd({ d,w_i) } ) }, u_j(\gcd({ d,w_j) }\gcd(d,w_i)
  ) \bigr) } {w_iw_j}}  \vspace{1mm} \\
 & = & {\displaystyle \frac{\gcd{ \bigl( d \gcd({d,w_j)}, 
d\gcd{(d,w_i)\bigr)} }}
{w_iw_j}=\frac{d\gcd{\bigl(\gcd{(d,w_j)},\gcd{ (d,w_i )\bigr)} }}
{w_iw_j} } \vspace{1mm} \\
 & = & {\displaystyle \frac{d\gcd{( w_i,w_j)  }}{w_iw_j} }\, ,
\end{array}
$$ 
where the last equality holds because
$\gcd{(w_i,w_j) }$ divides $d$. Alternatively (an argument that is
simpler), we could get the same result by showing that if 
$\gcd{( w_i,w_j ) } |d$, then
$${\rm lcm} \left( \frac{d}{{\gcd{( d,w_i) } }}, \frac{d}{
\gcd{( d,w_i ) } } \right) =\frac{d}{\gcd{( w_i,w_j) } } \, .$$ 

The last term is handled similarly,
$$\frac{u_1u_2u_3}{v_1v_2v_3{\rm
lcm}( u_1,u_2,u_3 ) }=\frac{d^2\gcd{( w_1,w_2,w_3 )} }{w_1w_2w_3}\, ,$$
because under the assumption that $\gcd{( w_i,w_j) }|d$ for all
$i\not=j$,  we have that 
$${\rm lcm}\left( \frac{d}{\gcd{ ( d,w_i) }},\frac{d}{\gcd{ ( d,w_i) }},
\frac{d}{\gcd{( d,w_k) }}\right) =\frac{d}{\gcd{( w_i,w_j,w_k) }}\, .$$ 
Combining all of these terms, we get that
$$2g_0=-1+\sum_{i=1}^3 \frac{\gcd{( d,w_i) }}{w_i}-
d\sum_{1\leq i<j\leq3}
\frac{\gcd{( w_i,w_j) }}{w_iw_j}+\frac{d^2\gcd{( w_1,w_2,w_3) }}
{w_1w_2w_3}\, .$$ 
If we now re-scale
$$\tilde{\bf w}=\frac{(w_1,w_2,w_3)}{\gcd{( w_1,w_2,w_3) }}\, ,\quad
\tilde{d}=\frac{d}{\gcd{( w_1,w_2,w_3) }}\, ,$$ 
we see that for $0\leq i<j<k\leq3$ we have that
$k_{ijk}=2g(\Sigma_{(\tilde{\bf w}, \tilde{d})})$, as desired.

It remains to show that $m_0=\gcd{( w_1,w_2,w_3) }$, which we
shall leave as an exercise.

The conjecture follows now by Proposition
\ref{prop.kol.links}. This proposition shows that when any of the
$m_i>1$, then the corresponding $g_i$ is the genus of a quasi-smooth
curve $D_i$, and therefore, $2g_i$ must not only be integral, but
also equal to twice that genus as it was shown. When $m_i=1$, it does
not matter what $2g_i$ is as it does not enter into the torsion
formula in either Koll\'ar's theorem or in Orlik's Algorithm.


\begin{thebibliography}{GM}
\bibitem{BM93}
A. Bangaya \& P. Molino, {\it G\'eom\'etric des formes de contact
compl\`etement int\'egrables de type toriques}, S\'eminaire G. Darboux
de G\'eom\'etrie et Topologie Diff\'erentielle, 1991-1992, Montpellier.
\bibitem{Bar65}
D. Barden, {\it Simply connected five-manifolds}, Ann. of Math. (2)
82 (1965), pp. 365-385. MR 32 \#1714.
\bibitem{BG00b}
C.P. Boyer \& K. Galicki, {\it A note on toric contact geometry}, J. Geom.
Phys. 35 (2000), 4, pp. 288--298. MR 2001h:53124.
\bibitem{BG06b}
\bysame, {\it Einstein metrics on rational homology
spheres}, J. Diff. Geom. (3) 74 (2006), pp. 353-362. MR MR2269781.
\bibitem{BG01b}
\bysame, {\it New Einstein metrics in dimension five}, J.
Diff. Geom. (3) 57 (2001), pp. 443-463. MR 2003b:53047.
\bibitem{BG00a}
\bysame, {\it On Sasakian-Einstein geometry},
Internat. J. Math. (7) 11 (2000),
pp. 873-909. MR 2001k:53081.
\bibitem{BG05}
\bysame, {\it Sasakian Geometry}, Oxford Mathematical
Monographs, Oxford University Press, to appear, 2008.
\bibitem{BoGa05a}
\bysame, {\it Sasakian geometry, hypersurface singularities, and Einstein
metrics}, Rend. Circ. Mat. Palermo (2) Suppl. (2005), 75, suppl., pp. 57-87.
 MR 2152356.
\bibitem{BGK05}
C.P. Boyer, K. Galicki \& J. Koll\'ar, {\it Einstein metrics on spheres},
 Ann. Math. (2) 162 (2005), no. 1, pp. 557-580. MR 2178969(2006j:53058).
\bibitem{BGKT05}
C.P. Boyer, K. Galicki, J. Koll{\'a}r \& E. Thomas, {\it Einstein metrics
on exotic spheres in dimensions 7, 11, and 15}, Experiment. Math. 14
(2005), no. 1, pp. 59-64. MR 2146519 (2006a:53042).
\bibitem{BGM06}
C.P. Boyer, K. Galicki \& P. Matzeu, {\it On eta-Einstein Sasakian geometry},
Comm. Math. Phys. (1) 262 (2006), pp. 177-208. MR 2200887.
\bibitem{BGO06}
C.P. Boyer, K. Galicki \& L. Ornea, {\it Constructions in Sasakian
Geometry}, preprint, arXiv:math.DG/0602233, to appear in Mathematische
Zeitscrift (2007).
\bibitem{bgs}
C.P. Boyer, K. Galicki \& S.R. Simanca, {\it Canonical Sasakian metrics},
Comm. Math. Phys (to appear), math.DG/0604325.
\bibitem{ca2}
E. Calabi, {\it Extremal K\"ahler metrics II}, in
Differential geometry and complex analysis (I. Chavel \&
H.M. Farkas eds.), Springer-Verlag, 1985, pp. 95-114.
\bibitem{CFO07}
K. Cho, A. Futaki \& H. Ono, {\it Uniqueness and examples of compact toric
Sasaki-Einstein metrics}, preprint; arXiv:math.DG/0701122 (2007).
\bibitem{fu}
A. Futaki, {\it An obstruction to the existence of Einstein K\"ahler
metrics}, Invent. Math., 73 (1983), pp. 437-443. MR 84j:53072. 
\bibitem{fow}
A. Futaki, H. Ono \& G. Wang,
{\it Transverse K\"ahler geometry of Sasaki manifolds and toric
Sasaki-Einstein manifolds}, preprint arXiv:math.DG/0607586, (2006).
\bibitem{GMSW04a}
J.P. Gauntlett, D. Martelli, J. Sparks \& W. Waldram, {\it Sasaki-Einstein
metrics on $S^2\times S^3$}, Adv. Theor. Math. Phys., 8 (2004), pp. 711-734.
\bibitem{GMSY06}
J.P. Gauntlett, D. Martelli, J. Sparks \& S.-T. Yau, {\it Obstructions to
the existence of Sasaki-Einstein metrics}, Comm. Math. Phys. (3) 273
(2007), pp. 803-827.
\bibitem{GhKo05}
A. Ghigi \& J. Koll\'ar, {\it K\"ahler-Einstein metrics on orbifolds and
Einstein metrics on Spheres}, Comment. Math. Helvetici (to appear),
arXiv:math.DG/0507289.
\bibitem{GrOr01}
G. Grantcharov \& L. Ornea, {\it Reduction of Sasakian manifolds},
J. Math. Phys. (8) 42 (2001), pp. 3809-3816. MR 1845220 (2002e:53060).
\bibitem{Fle00}
A.R. Iano-Fletcher, {\it Working with weighted complete intersections},
Explicit birational geometry of $3$-folds,
London Math. Soc. Lecture Note Ser. 281, pp. 101-173,
Cambridge Univ. Press, Cambridge, (2000).
\bibitem{JoKo01a}
J.M. Johnson \& J. Koll{\'a}r, {\it K\"ahler-Einstein metrics on log del
Pezzo surfaces in weighted projective 3-spaces}, Ann. Inst. Fourier
  (Grenoble) (1) 51 (2001), pp. 69-79. MR 2002b:32041.
\bibitem{Kol06a}
J. Koll\'ar, {\it Circle actions on simply connected 5-manifolds}, Topology
 (3) 45 (2006), pp. 643-671. MR 2218760.
\bibitem{Kol05b}
\bysame, {\it Einstein metrics on five-dimensional Seifert bundles},
J. Geom. Anal. (3) 15 (2005), pp. 445-476. MR 2190241.
\bibitem{Kol06}
\bysame, {\it Positive Sasakian structures on $5$-manifolds},
these Proceedings, Eds. Galicki \& Simanca, Birkhauser.
\bibitem{ls}
C. LeBrun \& S.R. Simanca, {\it On the K\"{a}hler Classes of Extremal
Metrics}, Geometry and Global Analysis (Sendai, Japan 1993),
First Math. Soc. Japan Intern. Res. Inst. Eds.
Kotake, Nishikawa \& Schoen.
\bibitem{Ler02a}
E. Lerman, {\it Contact toric manifolds}, J. Symp. Geom. (4) 1 (2002),
pp. 785-828, MR 2 039 164.
\bibitem{MaSpYau06}
D. Martelli, J. Sparks \& S.-T. Yau, {\it Sasaki-Einstein manifolds and
volume minimisation}, arXiv:hep-th/0603021, preprint 2006.
\bibitem{Mil68}
J.W. Milnor, {\it Singular points of complex hypersurfaces}, Annals of
Mathematics Studies 61, Princeton University Press, Princeton,
N.J., 1968. MR 39 \#969.
\bibitem{MiOr70}
J. Milnor \& P. Orlik, {\it Isolated singularities defined by weighted
homogeneous polynomials}, Topology 9 (1970), pp. 385-393. MR 45 \#2757.
\bibitem{Or72}
P. Orlik, {\it On the homology of weighted homogeneous manifolds}, Proceedings
of the $2^{\rm nd}$ Conference on Compact Transformation Groups (Univ. of
Mass., Amherst, Mass., 1971), Part I,
Lect. Notes Math. 298, Springer-Verlag (Berlin), 1972,
pp. 260-269. MR 55 \#3312.
\bibitem{Or72b}
\bysame, {\it Seifert manifolds}, Lect. Notes Math. 291, Springer-Verlag
(Berlin), 1972, MR 54 \#13950.
\bibitem{OrRa77a}
P. Orlik \& R.C. Randell, {\it The monodromy of weighted homogeneous
singularities}, Invent. Math. (3) 39 (1977), pp. 199-211. MR 57 \#314.
\bibitem{OrWa71a}
P. Orlik \& P. Wagreich, {\it Isolated singularities of algebraic surfaces
with $C^*$ action}, Ann. of Math. (2) 93 (1971), pp. 205-228. MR 44 \#1662.
\bibitem{Ran75}
R.C. Randell, {\it The homology of generalized Brieskorn manifolds},
Topology (4) 14 (1975), pp. 347-355. MR 54 \#1270.
\bibitem{Ruk95a}
P. Rukimbira, {\it Chern-Hamilton's conjecture and $K$-contactness},
Houston J. Math. (4) 21 (1995), pp. 709-718. MR 96m:53032.
\bibitem{si}
S.R. Simanca, {\it Canonical metrics on compact almost complex manifolds},
Publica\c{c}\~{o}es Matem\'aticas do IMPA, IMPA, Rio de Janeiro, 2004.
97 pp.
\bibitem{si2}
\bysame, {\it Heat Flows for Extremal K\"ahler Metrics}, Ann. Scuola Norm.
Sup. Pisa CL. Sci., 4 (2005), pp. 187-217.
\bibitem{si3}
\bysame, {\it Precompactness of the Calabi Energy}, Internat.
J. Math., 7 (1996) pp. 245-254.
\bibitem{si4}
\bysame, {\it Strongly Extremal K\"ahler Metrics}, Ann. Global Anal. Geom.
18 (2000), no. 1, pp. 29-46.
\bibitem{lusa1}
S.R. Simanca \& L.D. Stelling, {\it Canonical K\"ahler classes}. Asian J. Math.
5 (2001), no. 4, pp. 585-598.
\bibitem{Sm62}
S. Smale, {\it On the structure of $5$-manifolds}, Ann. of Math. (2)
  75 (1962), pp. 38-46. MR 25 \#4544.
\bibitem{Tak}
T. Takahashi, {\it Deformations of Sasakian structures and its application
to the Brieskorn manifolds}, T\^ohoku Math. J. (2) 30 (1978), no. 1,
pp. 37-43. MR 81e:53024.
\bibitem{WaZi90}
M.Y. Wang \& W. Ziller, {\it Einstein metrics on principal torus bundles},
J. Diff. Geom. (1) 31 (1990), pp. 215-248. MR 91f:53041.
\bibitem{YaYu05}
S.S.-T. Yau and Y. Yu, {\it Classification of $3$-dimensional isolated
rational hypersurface singularities with ${\mathbb C}^\ast$-action}, Rocky
Mountain J. Math. (5) 35 (2005), pp. 1795-1809. MR 2206037(2006j:32034).
\end{thebibliography}
\end{document}